\documentclass{article}
\usepackage[pdftex]{graphicx}
\pdfoutput=1
\usepackage{arxiv}
\usepackage[utf8]{inputenc} 
\usepackage[T1]{fontenc}    
\usepackage{hyperref}       
\usepackage{url}            
\usepackage{booktabs}       
\usepackage{amsfonts}       
\usepackage{nicefrac}       
\usepackage{microtype}      

\usepackage{doi}

\usepackage{multirow}
\usepackage{amsmath}
\usepackage{amssymb,xcolor}
\usepackage{caption}
\usepackage[linesnumbered,ruled,vlined]{algorithm2e}
\usepackage{algpseudocode}
\usepackage{enumerate}
\usepackage{xfrac}
\usepackage{amsthm}
\usepackage{comment}
\usepackage{subcaption}
\def\grad{\nabla}

\def\bc{\mathbf{c}}

\def\bx{\mathbf{x}}  

\def\bI{\mathbf{I}}

\def\L{{\boldsymbol{\Lambda}}}

\def\cA{\mathcal{A}}

\def\cD{\mathcal{D}}
\def\cE{\mathcal{E}}

\def\cG{\mathcal{G}}

\def\cL{\mathcal{L}}

\def\cN{\mathcal{N}}
\def\cO{\mathcal{O}}
\def\cP{\mathcal{P}}

\def\cS{\mathcal{S}}

\def\mP{\mathbb{P}}

\def\smskip{\smallskip}

\def\texitem#1{\par\smskip\noindent\hangindent 25pt
               \hbox to 25pt {\hss #1 ~}\ignorespaces}

\def\norm#1{\left\|#1\right\|}

\newcommand{\BEAS}{\begin{eqnarray*}}
\newcommand{\EEAS}{\end{eqnarray*}}
\newcommand{\BEA}{\begin{eqnarray}}
\newcommand{\EEA}{\end{eqnarray}}
\newcommand{\BEQ}{\begin{eqnarray}}
\newcommand{\EEQ}{\end{eqnarray}}
\newcommand{\BIT}{\begin{itemize}}
\newcommand{\EIT}{\end{itemize}}
\newcommand{\BNUM}{\begin{enumerate}}
\newcommand{\ENUM}{\end{enumerate}}

\newcommand{\BA}{\begin{array}}
\newcommand{\EA}{\end{array}}

\newcommand{\ones}{\mathbf 1}

\newcommand{\reals}{\mathbb{R}}
\newcommand{\integers}{\mathbb{Z}}

\newcommand{\diag}{\mathop{\bf diag}}

\newif\ifpagenumbering
\pagenumberingtrue

\pagenumberingfalse

\newsavebox{\theorembox}
\newsavebox{\lemmabox}
\newsavebox{\defnbox}
\newsavebox{\corollarybox}
\newsavebox{\propositionbox}
\newsavebox{\remarkbox}
\newsavebox{\assbox}
\savebox{\theorembox}{\noindent\bf Theorem}
\savebox{\lemmabox}{\noindent\bf Lemma}
\savebox{\defnbox}{\noindent\bf Definition}
\savebox{\corollarybox}{\noindent\bf Corollary}
\savebox{\propositionbox}{\noindent\bf Proposition}
\savebox{\remarkbox}{\noindent\bf Remark}
\savebox{\assbox}{\noindent\bf Assumption}

\newtheorem{theorem}{Theorem}  
\newtheorem{defi}[theorem]{Definition}
\newtheorem{proposition}[theorem]{Proposition}
\newtheorem{lemma}[theorem]{Lemma}

\newtheorem{remark}[theorem]{Remark}
\newtheorem{corollary}[theorem]{Corollary}
\newtheorem{result}{Result}

\newcommand{\beq}{\begin{equation}}
\newcommand{\eeq}{\end{equation}}
\newcommand{\beqa}{\begin{eqnarray}}
\newcommand{\eeqa}{\end{eqnarray}}
\newcommand{\beqas}{\begin{eqnarray*}}
\newcommand{\eeqas}{\end{eqnarray*}}

\newcommand{\hij}{H_{W}(i\to j)}
\newcommand{\htmij}{H_{\widetilde{M}}(i\to j)}
\newcommand{\hwprij}{H_{\overline{W}_{P^r}}(i\to j)}

\newcommand{\htw}{H_{W}(i\to j)}

\newcommand\str{\bgroup\markoverwith
{\textcolor{red}{\rule[0.5ex]{2pt}{1.5pt}}}\ULon} 
\def\mg#1{{#1}}

\def\bc#1{{#1}}

\def\tn{{\tilde{n}}}

\usepackage{soul}
\newcommand{\response}[1]{{#1}}

\def\id{\mathbf{I}}

\def\one{\mathbf{1}}

\def\sym{\mathbb{S}}
\def\sa#1{{#1}}

\def\L{{\cal L}}

\title{Randomized Gossiping with Effective Resistance Weights: Performance Guarantees and Applications}

\author{Bugra Can$^{*}$\\
        Management Sciences and Information Systems \\
        Rutgers Business School \\
        \texttt{bugra.can@rutgers.edu}
        \And 
        Saeed Soori\\
        Department of Computer Sciences\\
        University of Toronto \\
       \texttt{saeed.soori.sh@gmail.com}
        \And 
        Necdet Serhat Aybat\\
        Industrial and Manufacturing Engineering Department\\
        Penn State University\\
        \texttt{nsa10@psu.edu}
        \And
        Maryam Mehri Dehnavi\\ 
        Department of Computer Sciences\\
        University of Toronto \\
        \texttt{mmheride@cs.toronto.edu}
        \And
        Mert G\"urb\"uzbalaban \thanks{Bugra Can and Mert Gürbüzbalaban acknowledge support from the Office of Naval Research Award Number N00014-21-1-2244, and the grants National Science Foundation (NSF) CCF-1814888, NSF DMS-2053485, NSF DMS-1723085.} \\
        Management Sciences and Information Systems \\
        Rutgers Business School \\
        \texttt{mert.gurbuzbalaban@rutgers.edu}
}

\renewcommand{\shorttitle}{Randomized Gossiping with Effective Resistance Weights}

\hypersetup{
pdftitle={Randomized Gossiping with Effective Resistance Weights: Performance Guarantees and Applications},
pdfsubject={},
pdfauthor={Bugra Can, Saeed Soori, Necdet Serhat Aybat, Maryam Mehri Dehnavi, Mert G\"urb\"uzbalaban},
pdfkeywords={Distributed algorithms/control, networks of autonomous agents,optimization, randomized gossiping algorithms},
}

\begin{document}
\maketitle

\begin{abstract}
	The effective resistance between a pair of nodes {in} a weighted undirected graph is defined as the potential difference induced 
when a unit current is injected at one node and extract from the other, {treating edge weights as the conductance values of edges}. The effective resistance is a key quantity of interest in many applications, e.g., solving linear systems, Markov Chains, and continuous-time averaging networks. We consider effective resistances~(ER) in the context of designing randomized gossiping methods for the consensus problem, where the aim is to compute the average of node values in a distributed manner through iteratively computing weighted averages among randomly chosen neighbours. 
\mg{For barbell graphs}, we \response{prove} that \response{choosing wake-up and communication probabilities proportional to }
ER weights improves the averaging time corresponding to the traditional choice of uniform weights. \mg{For $c$-barbell graphs, we show that ER weights admit lower and upper bounds on the averaging time that improves upon the lower and upper bounds available for uniform weights. \sa{Furthermore,} for graphs with a small diameter, we can show that ER weights can improve upon the existing bounds for Metropolis weights by a constant factor under some assumptions.} We illustrate these results through numerical experiments \mg{where we showcase the efficiency of our approach on several graph topologies including barbell graphs, small-world graphs, and stochastic block models}. We also present an application of the ER gossiping to distributed optimization: we numerically \response{verify} that using ER gossiping within EXTRA and DPGA-W methods improves their practical performance in terms of communication efficiency.
\end{abstract}

\keywords{Distributed algorithms/control \and networks of autonomous agents \and optimization \and randomized gossiping algorithms}

\section{Introduction}
\label{sec:intro}
\noindent
Let $\cG=(\mathcal{N},\mathcal{E},w)$ be an undirected, weighted and connected graph defined by the set of nodes (agents) $\cN=\{1,\ldots,n\}$, the set of edges $\cE\subseteq \cN\times\cN$, and 
the edge weights $w_{ij} > 0$ for $(i,j) \in \cE $. Since $\cG$ is undirected, we assume that both $(i, j)$ and $(j, i)$ refer to the same edge when it exists, and for all $(i,j)\in\mathcal{E}$, we set $w_{ji} = w_{ij}$. Identifying the weighted graph $\cG$ as an electrical network 
in which each edge $(i,j)$ corresponds to a branch of conductance $w_{ij}$, the effective resistance $R_{ij}$ between a pair of nodes $i$ and $j$ is defined as the voltage potential difference induced between them when a unit current is injected at $i$ and extracted at $j$.
The effective resistance (ER), also known as the resistance distance, is a key quantity of interest to compute in many applications and algorithmic questions over graphs. It defines a metric on \response{the graph} providing bounds on its conductance \cite{klein2002resistance,Klein1993}. Furthermore, it is closely associated with the hitting and commute times 
\bc{of} a random walk\footnote{The hitting time is the
expected number of steps of a random walk starting from $i$ until it first visits $j$. The commute time $C_{ij}$ is the expected number of steps required to go from $i$ to $j$ and from $j$ to $i$ back again.} on the graph $\mathcal{G}$ 
when the probability of a transition from $i$ to $j\in\cN_i$ is $w_{ij}/\sum_{j'\in\cN_i}w_{ij'}$ where {$\cN_i\triangleq\{j\in\cN:~ w_{ij}>0\}$} denotes the set of neighboring nodes of $i\in\cN$; therefore, it arises naturally for studying random walks over graphs and their mixing time properties \cite{boydSIAMreview,aldous-fill-2014,doyle1984random}, spectral approximation of graphs \cite{SpielmanSrivastava}, continuous-time averaging networks including consensus problems in distributed optimization \cite{boydSIAMreview}. 

There exist \emph{centralized} algorithms for computing or approximating effective resistances accurately which require global communication beyond local 
information exchange among the neighboring agents \cite{Mishra_2020,Jafarizadeh_2008,Jafarizadeh2006CalculatingTR, SpielmanSrivastava,bapat2003simple}. \mg{The references \cite{Mishra_2020,Jafarizadeh_2008,Jafarizadeh2006CalculatingTR} develop key techniques for computing the effective resistances explicitly on specific network types. In particular, \cite{Jafarizadeh_2008} addresses a class of graphs which are underlying networks of some symmetric association schemes whereas \cite{Mishra_2020} considers two dimensional resistor networks. The reference \cite{Jafarizadeh2006CalculatingTR} provides an algorithm for the calculation of the resistance between two arbitrary nodes in a distance-regular network and also provides analytical formulas. The works \cite{SpielmanSrivastava,bapat2003simple}} are based on computing or approximating the entries of the pseudoinverse $\cL^{\dag}$ of the graph Laplacian matrix $\cL$, based on the identity~\cite{SpielmanSrivastava}
{\small
\begin{align}
    R_{ij} =\cL^{\dag}_{ii} +\cL^{\dag}_{jj} - 2\cL^{\dag}_{ij},\quad \forall~(i,j)\in\cE.
\end{align}}%
However, such centralized algorithms are \emph{impractical} or \emph{infeasible} for several key applications in multi-agent systems, e.g., randomized gossiping algorithms, for averaging the node values across the whole network, 
use only local communications between random neighbors 
(see \cite{nedic2009distributed,boyd2006randomized,olfati2007consensus}); this motivates the use of \emph{distributed algorithms} for computing effective resistances {which only rely on the information exchange among immediate neighbors.}
{In these applications,  communication 
among the agents is typically the bottleneck 
compared to the complexity of local computations 
of the agents; 
thus, it is crucial to develop distributed algorithms that are efficient in terms of \response{the} total number of communications required.} 
To the best of \response{the} authors' knowledge, the first attempt for computing effective resistances in a decentralized way and also the first 
ER-based {randomized} gossiping algorithms appeared in \cite{aybat2017decentralized}. \bc{The latter algorithms are asynchronous gossiping algorithms where each agents' wake-up and communication probabilities are chosen proportional to ER weights (see Section \ref{sec:prelim} for details)}. Aybat and G\"urb\"uzbalaban 
have shown in \cite{aybat2017decentralized} that effective resistance~(ER) weights can be computed at each agent locally with an efficient distributed algorithm, Distributed Randomized Kaczmarz (D-RK). 
Our paper is motivated by the numerical evidence presented in~\cite{aybat2017decentralized} that using ER weights has the potential to improve the performance of randomized gossiping algorithms on specific graphs. Since in~\cite{aybat2017decentralized} no rigorous performance guarantees for the use of ER weights \response{were} provided, here we focus on establishing the missing theoretical results that match the outstanding empirical behavior.

\indent \textbf{Contributions.}
First, in this paper, we provide theoretical guarantees on the 
ER-based randomized gossiping algorithms proposed in \cite{aybat2017decentralized} for the consensus problem, where the objective is to compute the average of node values over a network in a 
decentralized manner \cite{boyd2006randomized}. {A standard approach for solving the consensus problem is the \emph{randomized uniform gossiping}}, where each node keeps a local estimate of the average of node values and has the equal (uniform) probability of being activated to communicate with a randomly chosen neighbour to update its local estimate. However, this approach treats all the edges (equally) uniformly and can be slow in practice. To overcome this problem, in \cite{aybat2017decentralized}, 
ER-based randomized gossiping algorithms were proposed without any theoretical guarantees, in which the edges are being activated by non-uniform probabilities that are proportional to their effective resistances. 

Our theoretical results presented in Section~\ref{sec: Main Results} (see Results~\ref{res: c-barbell graph},~\ref{res: eig-order-imprv}, and~\ref{res: d_graph_results}) explain the superior empirical behaviour of ER-based gossiping \response{over the uniform gossiping} observed in~\cite{aybat2017decentralized}. 
Briefly, we bound the time required to compute an inexact average using \response{analysis} based on conductance and spectral properties of the underlying weighted communication graph, and compare the bounds we obtained corresponding to
\response{the} ER and uniform gossiping methods. 
We 
show that averaging time with 
ER weights is $\Theta(n)$ faster than that of 
uniform gossiping on a barbell graph where $n$ is the number of agents. Furthermore, we also prove that for 
connected graphs with a small diameter, the averaging time with resistance weights can be faster than known performance bounds for the averaging time with gossiping based on Metropolis weights \mg{by a constant factor} \mg{see (Remark \ref{remark-metropolis-compare})}. \bc{We also provide numerical experiments on several graph topologies which illustrate the performance improvements that can be obtained within ER-based gossiping. In our experiments, 
the effective resistances are first computed with the normalized D-RK algorithm of \cite{aybat2017decentralized} and then used for ER-based gossiping.} \mg{Our theoretical and numerical results show that ER weights are especially useful in the presence of ``bottleneck edges" or \sa{clusters 
giving a graph cut leading to small graph conductance values.}}

On a different note, Aybat and G\"urb\"uzbalaban~\cite{aybat2017decentralized} introduced two alternative methods to compute ER weights in a decentralized manner: 
D-RK and normalized D-RK --both converging linearly.
In our experiments at Section~\ref{sec-numerical}, we 
have adopted the normalized D-RK, upon proving that the convergence rate of normalized D-RK is better than D-RK; {resolving a conjecture raised 
in~\cite{aybat2017decentralized}} (see the Supplementary Material). 

Second, we consider the consensus optimization problem, where the agents 
connected on a network aim to collaboratively solve the optimization problem {$\min_{x\in\mathbb{R}^p} f(x) \triangleq\sum_{i=1}^n f_i(x)$} where $f_i(x):\mathbb{R}^p \to \mathbb{R}$ is a cost function only available to (node) agent $i$. This problem includes a number of key problems in supervised learning including distributed regression and logistic regression or more generally distributed empirical risk minimization problems \cite{EmpRisk_Xiao,EmpRisk_Lee}. The consensus iterations are a building block of many existing state-of-the-art distributed consensus optimization algorithms such as the EXTRA 
and the distributed proximal gradient 
(DPGA-W) \cite{aybat2018distributed} algorithms for consensus optimization. We show through numerical experiments that our 
framework based on effective resistances can {improve the performance of} 
the EXTRA and DPGA-W algorithms for consensus optimization {in terms of the total number of communications required}. We believe our framework has far-reaching potential for improving the communication efficiency of many other distributed algorithms including distributed subgradient and ADMM methods, and this will be the subject of future work.

\textbf{Related work.} For consensus problems, there are some alternative methods to accelerate the commonly used consensus protocols. The approach in \cite{olshevsky2014linear} is a synchronous algorithm combining Metropolis weights with a momentum averaging scheme. There are other approaches based on momentum averaging \cite{loizou2018accelerated,loizouRevisiting,loizouAccConsensus}, min-sum splitting \cite{RebeschiniMinSumSplitting}, and Chebyshev acceleration \cite{ConsensusAccelerateChebyshev,SeidmanChebyshev,bu2018accelerated} to accelerate the convergence speed of the consensus methods. 
This paper is orthogonal to the \mg{momentum averaging-based} approaches \response{in the sense that} it can be used in combination with \response{the} aforementioned \mg{momentum-based} acceleration schemes, \mg{we refer the reader to the Supplementary Material for the details. There are also works that provide lower bounds on the distributed averaging time on a graph \cite{boyd2006randomized,Bound_Roch,Boyd03fastestmixing,shah-book}. In particular, it follows from these lower bounds that for the two-dimensional grid, even the best gossiping weights will not lead to an accelerated performance compared to baseline approaches. Indeed, for special graphs such as the two-dimensional grid, cycle graph or the line graph, ER weights will be similar to uniform weights due to the symmetries in the graph structure and consequently ER weights will not improve the performance compared to uniform weights. However, for graphs with asymmetries involving clusters or bottleneck edges along which the graph cut has low conductance, based on our numerical and theoretical results, we expect ER weights lead to an improved performance.}  
\\
\indent \textbf{Outline.} 
In Section~\ref{sec:prelim},
we give a brief overview of randomized gossiping including uniform and ER-based gossiping methods.
In Section~\ref{sec: Main Results}, we state our main contributions.
In Section~\ref{sec-theoretical guarantees}, we provide detailed arguments establishing the main results stated in Section~\ref{sec: Main Results}. 
In Section \ref{sec-numerical}, we provide numerical experiments 
illustrating that 
using 
ER weights can improve the performance of EXTRA and DPGA-W algorithms for consensus optimization. 
In Section~\ref{sec-future}, we give some concluding remarks. Finally, we present some of the proofs and supporting results in Appendix \ref{sec: Prop 6and9}--\ref{sec: Supporting Results}.

{\textbf{Notation.} Let $|S|$ denote the cardinality of a set $S$, $\lfloor . \rfloor$ denote the floor function and $\integers_+$ be the set of nonnegative integers}. We define $d_i\triangleq |\cN_i|$ as the degree of $i\in\cN$, and $m\triangleq|\mathcal{E}|$.  Throughout the paper, $\cL\in\reals^{|\cN|\times |\cN|}$ denotes the weighted Laplacian of $\cG$, i.e., $\L_{ii}=\sum_{j\in\cN_i}w_{ij}$, $\cL_{ij}=-w_{ij}$ if $j\in\cN_i$, and equals to $0$ otherwise. The diameter of a graph is $\mathcal{D}\triangleq\max_{i,j\in \mathcal{N}}d(i,j)$ where $d(i,j)$ is the shortest path on the graph between nodes $i$ and $j$. 
The set $\sym^n$ denotes the set of $n\times n$ real symmetric matrices. 
We use the notation $Z = [z_i]_{i=1}^n$ where $z_i$'s are either the columns or rows of the matrix $Z$ depending on the context. 
$\ones$ is the column vector with all entries equal to 1, and $\id$ is the identity matrix. We let $||x||_p$ denote the $L_p$ norm of a vector $x$ for $p\geq 1$, and let $\|A \|_F$ denote the Frobenius norm of a matrix $A$. A square matrix $A$ is \emph{doubly stochastic} if all of its entries are non-negative and all its rows and columns sum up to 1. We say that a square matrix $A$ is \emph{weakly diagonally dominant} if it's diagonal entries $A_{ii}$ satisfy the inequality $|A_{ii}|\geq \sum_{j\neq i}|A_{ij}|$ for every $i$. Let $f$ and $g$ be real-valued functions defined over positive integers. We say $f(n)=\cO(g(n))$ if $f$ is bounded above by $g$ asymptotically, i.e., there exist constants $ k_{1}>0$ and $n_{0}\in\integers_+$ such that 
$ f(n)\leq k_{1}\cdot g(n)$ for all $n>n_{0}$. Similarly, we say $f(n)=\Omega (g(n))$ if there exist constants $k_2 > 0$ and $n_0\in\integers_+$ such that $f(n)\geq k_2 \,g(n)$ for every $n > n_0$; and we say $f(n)=\Theta(g(n))$ if $f(n)=\Omega (g(n))$ and $f(n)=\cO (g(n))$.
Finally, $\log(x)$ denots the natural logarithm of $x$, and $e_i$ is the $i$-th standard basis vector in $\mathbb{R}^n$ for $i=1,2,\dots,n$. 
\vspace{-3mm}
\section{Preliminaries}
\label{sec:prelim}
\subsection{Randomized gossiping}
\label{sec:gossiping}
Here we 
give an overview of randomized gossiping methods 
for the consensus problem. These methods can compute the average of node values over a network in an asynchronous and decentralized manner, for 
details 
see \cite{boyd2006randomized,shah-book}.

Let $y^0 \in \mathbb{R}^n$ be a vector such that the $i$-th component {$y_i^0$} 
represents the initial value at node $i\in\cN$. 
The aim of the randomized gossiping algorithms is to have each node compute the average $\bar{y} \triangleq \sum_{i=1}^n y_i^0/n$ in a 
decentralized manner through an iterative procedure. 
At every iteration $k\in\integers_+$, each node $i\in\cN$ possesses a local estimate $y_i^k$ of the average to be computed and communicates with only randomly selected neighbors to update its estimate.  
{The setup is that} 
each node $i\in\cN$ has an exponential clock ticking with rate $r_i>0$ where 
the time between two ticks is exponentially distributed and independent 
of other nodes' clocks. 
A node wakes up when its clock ticks. \bc{Since all the clocks are independent,} if a node wakes up at time $t_k\geq 0$, it is node $i$ with probability (w.p.) $p_i\triangleq r_i/\sum_{j\in\cN}r_j$. Given that 
the node $i$ wakes up at time $t_k$, the conditional probability that it picks \emph{one} of its neighbors $j\in\cN_i$ to communicate with probability $p_{j|i}\in(0,1)$, where the probabilities $\{p_{j|i}\}_{j\in\cN_i}$ are design parameters satisfying $\sum_{j\in\cN_i}p_{j|i}=1$. When either $i$ wakes up and picks $j\in\cN_i$ or vice versa, we say the edge $(i,j)$ is activated. Once the edge $(i,j)$ is activated, nodes $i$ and $j$ exchange their local variables $y_i^k$ and $y_j^k$ at time $t_k$ and both compute the average $(y_i^k+y_j^k)/2$. {This is illustrated in Algorithm \ref{alg:gossip} below which admits an asynchronous implementation -- see, e.g., \cite{boyd2006randomized}}. 
\DontPrintSemicolon
\begin{algorithm}
\small
\textbf{Initialization:} $y^{0}=[y_1^{0}, y_2^{0}, \dots, y_n^{0}]^\top \in\mathbb{R}^n$ \;
\For{$k\geq 0$}{
At time $t_k$, $i\in\cN$ wakes up w.p. $p_i=r_i/\sum_{j\in\cN}r_j$\;
Picks $j\in\cN_i$ randomly w.p. $p_{j|i}$\;
$y_i^{k+1}\gets \frac{y_i^{k}+y_j^{k}}{2}$,\quad $y_j^{k+1}\gets \frac{y_i^{k}+y_j^{k}}{2}$
\label{eq:y}\;
}
\caption{\small Randomized Gossiping}
\label{alg:gossip}
\end{algorithm}
{ Assuming {\it there are no self-loops }
for each $i\in\cN$, let
{\small
\begin{subequations}
\label{eq:P}
\begin{align}
    &P_{ii}\triangleq 0;\quad P_{ij} \triangleq p_i~p_{j|i},\quad \forall\ j\in\cN_i,\\ 
    &P_{ij}\triangleq 0,\quad \forall j\in\cN\setminus\cN_i,
\end{align}
\end{subequations}}%
where $P_{ij}$ is the (unconditional) probability that the edge $(i,j)$ is activated by the node $i$. 
By definition, we have $\sum_{ij} P_{ij} \triangleq \sum_{i\in\cN}\sum_{j\in\cN} P_{ij}= 1$. Let $\mathcal{A}(P)$ denote an asynchronous gossiping algorithm characterized by a probability matrix $P$ as in~\eqref{eq:P} for some set of probabilities $\{p_i\}_{i\in\cN}$ and $\{p_{j|i}\}_{j\in\cN_i}$ for $i\in\cN$. The performance of 
$\cA(P)$ is typically measured by the $\varepsilon$-averaging time, defined for any $\varepsilon>0$ as}:
{\small
\begin{align}
\label{eq:T_ave}
    T_{ave}(\varepsilon,P)\triangleq \sup_{y^0\in \mathbb{R}^n\setminus\{\bf 0\}}\inf \left\{k : \mP\left( \frac{\|y^k - 
\bar{y}\one \|}{\|y^0\| } \geq \varepsilon \right) \leq \varepsilon \right\},
\end{align}}%
{see, e.g., \cite{boyd2006randomized}}. Suppose 
$(i,j)$ is 
activated by node $i$, then we can write the update in Step~\ref{eq:y} of the Algorithm~\ref{alg:gossip} as
{\small
$$ y^{k+1} = \bc{W^{(i,j)}} y^k \quad
\mbox{where} \quad \bc{W^{(i,j)}} \triangleq I - \frac{(e_i - e_j)(e_i - e_j)^\top}{2}. $$}%
We also define
    \begin{equation}\label{def-Wp}
        \overline{W}_P \triangleq \mathbb{E}_P [\bc{W^{(i,j)}}] = \sum_{i,j\in \mathcal{N}} P_{ij} \bc{W^{(i,j)}},
    \end{equation}
which is the expected value of the random iteration matrix \response{$W^{(i,j)}$} with respect to the distribution defined over $i\in\cN$ and  $j\in\cN_i$. The following theorem from \cite{boyd2006randomized} shows that the second largest eigenvalue of $\overline{W}_P$ determines the $\varepsilon$-averaging time.
\begin{theorem}[{\cite[Theorem 3]{boyd2006randomized}}]\label{thm-gossip-compare} For a given $\cA(P)$, the symmetric matrix $\overline{W}_P$ defined in \eqref{def-Wp} satisfies 
{\small
  \begin{align*} 
      0.5 \frac{\log(\varepsilon^{-1})}{\log ([\lambda_{n-1}(\overline{W}_P)]^{-1})}  &\leq T_{ave}(\varepsilon,P) \leq 3 \frac{\log(\varepsilon^{-1})}{\log ([\lambda_{n-1}(\overline{W}_P)]^{-1})},
  \end{align*}}%
where $\lambda_{n-1}(\overline{W}_P)$ is the second largest eigenvalue of $\overline{W}_P$.
\end{theorem}
This result makes the connection between the convergence time of an asynchronous gossiping algorithm $\mathcal{A}(P)$ and the spectrum of the expected iteration matrix $\overline{W}_P$. It is therefore of interest 
to design $P$ through carefully choosing the probabilities $\{p_i\}_{i\in\cN}$ and $\{p_{j|i}\}_{j\in\cN_i}$ for $i\in\cN$ in order to get the best performance, i.e., the smallest $\varepsilon$-averaging time. 

In this paper, we consider two different randomized gossiping algorithms: \emph{uniform gossiping} and \emph{ER gossiping} which differ in how the probabilities $\{p_i\}_{i\in\cN}$ and $\{p_{j|i}\}_{j\in\cN_i}$ for $i\in\cN$ are selected.  
In particular, based on Theorem \ref{thm-gossip-compare}, we will study the second largest eigenvalue of the expected iteration matrix $\overline{W}_P$ corresponding to these two algorithms and compare their $\epsilon$-averaging times.

\subsection{Randomized uniform gossiping}
\begin{figure}[h!]
\begin{center}
	\includegraphics[width=0.4\linewidth]{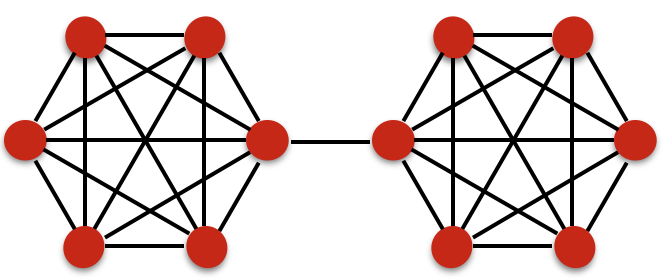}
    \caption{\label{fig-barbell}\small Barbell graph $K_{\tilde{n}} - K_{\tilde{n}}$ with $n=2\tilde{n}=12$ nodes}
\end{center}
\end{figure}

{{In the \textit{randomized uniform gossiping}, 
each node $i$ wakes up with equal probability $p_i^u=\frac{1}{n}$, i.e., using uniform clock rates $r_i=r>0$ for $i\in\cN$. The superscript $u$ stands for the uniform choice of clock rates. Then, node $i$ picks the edge $(i,j)$ with conditional probability $p_{j|i}^u=\frac{1}{d_i}$ for 
$j\in\cN_i$; thus, 
\begin{equation*}
    P^u_{ij}=p^u_ip^u_{j|i}=\frac{1}{n}\frac{1}{d_i},
\end{equation*}
see, e.g., \cite{afek2013distributed,xiao2004fast}}.}
{One of the drawbacks of this approach is that it can be quite slow over graphs with a high \emph{bottleneck ratio} \cite{fagnani2017introduction} where, intuitively speaking, some 
``bottleneck edges" 
limit the spread of information over the underlying graph. A classical example of a graph with a high bottleneck ratio is the \emph{barbell graph}. Barbell graphs are frequently studied 
within the consensus problem literature as they constitute a worst-case example in terms of both the mixing properties of random walks~\cite[Section 5]{aldous-fill-2014} and the performance of distributed averaging algorithms (see, e.g., \cite{boydSIAMreview,jung2010distributed}). 

Barbell graphs 
consist of two complete subgraphs connected 
with an edge (see Figure \ref{fig-barbell}). Let $K_\tn$ denote a complete graph with $\tn$ nodes, we will be denoting a barbell graph with $n=2\tn$ nodes by $K_\tn - K_\tn$. Let $(i^*,j^*)$ be the edge that connects the two complete subgraphs which we will be referring to as the \emph{bottleneck edge}. This is the only edge that allows node values to be propagated between the two complete subgraphs; therefore, how frequently it is sampled is a key factor that determines the averaging time.}

{
The 
probability of sampling the \emph{bottleneck edge} $(i^*,j^*)$, with uniform weights can be computed explicitly:

\BEQ \label{eq-bottleneck-uniform-proba} 
P^u_{i^*j^*} = P^u_{j^*i^*} = \frac{1}{n} \frac{1}{d_{i^*}} = \response{ \frac{2}{n^2}}.
\EEQ

}
{This implies that it takes $\Theta(n^2)$ iterations in expectation to activate this edge, which 
is the underlying reason why the randomized uniform gossiping 
iterates converge slowly when $n$ is large on the barbell graph}. {
{The effect of bottleneck edges on the performance of gossiping algorithms has been recently studied experimentally by Aybat and G\"urb\"uzbalaban \cite{aybat2017decentralized} on different 
topologies
including the barbell and small-world graphs. 
The authors proposed 
\emph{ER gossiping} where the edges are sampled with non-uniform probabilities proportional to effective resistances $\{R_{ij}\}_{(i,j)\in\cE}$ and the numerical experiments in~\cite{aybat2017decentralized} showed that this can lead to significant performance improvement over graphs 
with bottleneck edges, 
such as barbell graphs. We next describe this method.}}
\vspace*{-4mm}
\subsection{Effective-resistance~(ER) gossiping}

{
In the ER 
gossiping, 
each $i\in\cN$ wakes up with probability $p_i^r=\frac{\sum_{j\in\cN_i}R_{ij}}{2\sum_{(i,j)\in\cE}R_{ij}}$, i.e., 
setting clock rate $r_i=\sum_{j\in\cN_i}R_{ij}$ for $i\in\cN$, and node $i$ picks $(i,j)$ with conditional probability $p_{j|i}^r=\frac{R_{ij}}{\sum_{j\in\cN_i}R_{ij}}$ for all $j\in\cN_i$; thus, 
ER gossiping corresponds to the unconditional probabilities
$$P^r_{ij} = p_i^r p_{j|i}^r = \frac{R_{ij}}{2\sum_{(i,j)\in\cE} R_{ij}}= \frac{R_{ij}}{2(n-1)}=P_{ji}^r, 
$$
for all $(i,j)\in \cE$ where the third equality follows from Foster’s Theorem which says that $\sum_{(i,j)\in\cE} R_{ij}=(n-1)$ -- see, e.g., \cite{tetali1994extension}. This choice of sampling probabilities can lead to bottleneck edges being more frequently sampled. 
}
{We illustrate this fact on the barbell graph ($K_\tn-K_\tn$): Note that the unconditional probability of sampling the bottleneck edge $(i^*, j^*)$ is given explicitly as 
\BEQ\label{eq-bottleneck-resist-proba} P^r_{i^*j^*} = P^r_{j^*i^*} = \frac{R_{i^*j^*}}{2(n-1)}= \frac{1}{2(n-1)}, 
\EEQ
where $n=2\tn$ and we used the fact that $R_{i^* j^*}=1$ (see the proof of Lemma \ref{lemma: Values} for the derivation of \eqref{eq-bottleneck-resist-proba}). Hence, comparing \eqref{eq-bottleneck-uniform-proba} and \eqref{eq-bottleneck-resist-proba}, we see that ER weights allow sampling of the bottleneck edge $(i^*,j^*)$ more frequently, by a factor of $\Theta(n)$, than the uniform gossiping} on $K_\tn-K_\tn$. Intuitively speaking, this is the reason why ER 
gossiping can be efficient on barbell graphs. Numerical experiments provided {in~\cite{aybat2017decentralized} support this intuition where ER gossiping 
outperforms uniform gossiping over an unweighted barbell graph as well as small-world graphs, which are random graphs that arise frequently in real-world applications such as social networks. 
}

\begin{figure}[ht!]
    \centering
    \includegraphics[width=.7\linewidth]{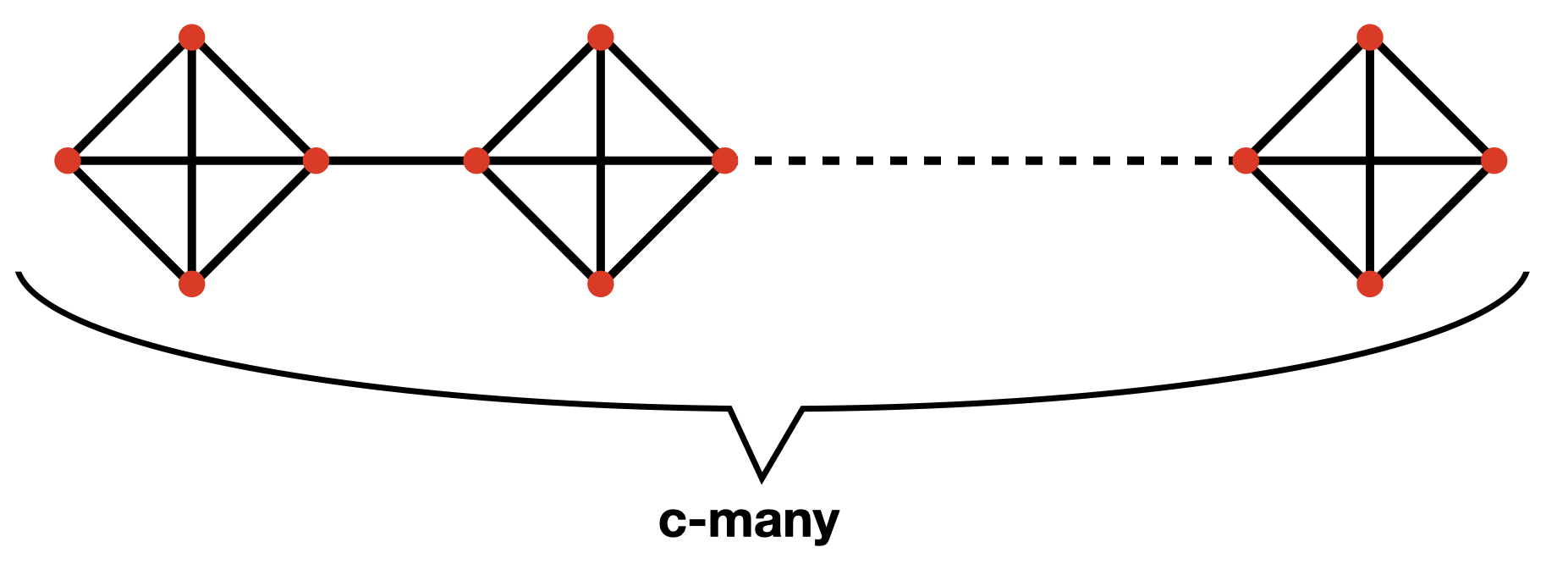}
    \caption{A $c$-barbell graph with $\tn=4$ (
    $n=\tilde{n} c = 4c$).}
    \label{fig: Example c-Kn}
    \vspace*{-4mm}
\end{figure}

Despite the empirical success of \emph{ER gossiping} in practice, theoretical results 
supporting its practical performance have been lacking in the literature. The purpose of this paper is to provide rigorous convergence guarantees for ER gossiping algorithms on certain network topologies (see Section~\ref{sec: Main Results} for our main results' statements and Section~\ref{sec-theoretical guarantees} for the proofs) and 
to present further numerical evidence that ER gossiping, beyond distributed averaging, can also improve the practical performance of distributed methods for consensus optimization (Section~\ref{sec-numerical}). Indeed, in our analysis, we consider connected graphs characterized by their diameter $\mathcal{D}\in \mathbb{Z}_+$, barbell graphs and $c$-barbell graphs which are generalizations of barbell graphs. More specifically,  a $c$-barbell graph ($ K_{\tn}^c$) for $c\geq 2$ is a path of $c$ equal-sized complete graphs ($K_{\tn}$) \cite{censor2012fast}, e.g., see 
Figure~\ref{fig: Example c-Kn} for $K_4^c$. In the special case, when $c=2$, a $c$-barbell graph is equivalent to the barbell graph. We show that for these graphs, ER  gossiping has provably better convergence properties than uniform gossiping in terms of $\varepsilon$-averaging times. Precise results will be stated in the next section.
\section{Main Results} \label{sec: Main Results}
{In this section, we state our main theoretical results: we provide performance bounds 
for the ER 
gossiping} in terms of  $\epsilon$-averaging time $T_{ave}(\varepsilon,P^r)$. {Our results highlight the performance improvements 
obtained with this approach.} 
{Our first result concerns $c$-barbell graphs where we focus on the $\varepsilon$-averaging times 
of uniform 
and ER 
gossiping algorithms}. To the best of our knowledge, for $c$-barbell graphs, an analytical formula for the second largest eigenvalue $\overline{W}_P$ is not analytically available; therefore, in our analysis we estimate this eigenvalue based on graph conductance techniques (see Section~\ref{Subsec: Conductance} for details) which leads to the following lower and upper bounds on the $\varepsilon$-averaging times.

{
\begin{result}\label{res: c-barbell graph} Given $\epsilon>0$, and $\tilde{n},c\in\integers_+$ such that $c\geq 2$, asynchronous randomized gossiping algorithms $\mathcal{A}(P^u)$ and $\mathcal{A}(P^r)$ on a c-barbell graph with $n=\tilde{n}c$  
satisfy 
\begin{align}
\Theta(c^2\tilde{n}^3\log(1/\epsilon)) \leq T_{ave}(\varepsilon,P^u ) \leq \Theta(c^4{\tilde{n}}^6\log(1/\epsilon)), \label{bound: c_barbell_uni}
\\
\Theta(c^2{\tn}^2\log(1/\epsilon)) \leq T_{ave}(\varepsilon,P^r ) \leq\Theta(c^4{\tn}^4\log(1/\epsilon)). \label{bound:c_barbell_res}
\end{align}
\end{result} 

These bounds from 
 Result~\ref{res: c-barbell graph} for the c-barbell graph show that, for any given precision $\epsilon>0$, using effective resistances 
  one can improve upper and lower bounds on the averaging times 
  by a factor of $\Theta(n)$ and $\Theta(n^2)$, respectively. \mg{In the Supplementary Material, we also compared the averaging times $T_{ave}(\varepsilon,P^r)$ and $T_{ave}(\varepsilon,P^u )$ numerically based on computing the second-largest eigenvalues $\lambda_{n-1}(\overline{W}_{P^r})$ and $\lambda_{n-1}(\overline{W}_{P^u})$ and by invoking Theorem \ref{thm-gossip-compare}. These numerical results are inline with Result \ref{res: c-barbell graph}, showing that effective resistances improve upon uniform weights in the sense that the averaging time for the effective resistances scales better with the number of nodes $\tilde n$.}
  
  The next result shows that for the case of barbell graphs (when $c=2$) 
the 
ER gossiping is in fact faster by a factor of $\Theta(n)$. The proof idea is based on computing the eigenvalues of 
$\overline{W}_{P^r}$ and $\overline{W}_{P^u}$ explicitly via exploiting symmetry group properties of barbell graphs and showing that the lower bounds in \eqref{bound: c_barbell_uni}--\eqref{bound:c_barbell_res} 
are 
{attained} for $c=2$.

\begin{result}\label{res: eig-order-imprv}Given $\epsilon>0$ and $n\in\integers_+$, let $n=2\tn$. The $\varepsilon$-averaging times of asynchronous gossiping algorithms $\mathcal{A}(P^r)$ and $\mathcal{A}(P^u)$ on barbell graph $K_{\tn}-K_{\tn}$ satisfy the equality:
    $$T_{ave}(\varepsilon,P^r) = \Theta(1/n)~ T_{ave}(\varepsilon, P^u).$$
\end{result}
{
A natural question is whether it is possible to further improve the ER 
gossiping bounds for barbell graphs; however, in the next result, we show that this is not possible as long as the matrix $P$ is symmetric --thus, ER gossiping is optimal. Finally, we also obtain $\varepsilon$-averaging bounds for a more general class of connected graphs 
depending on their diameters.} 

\begin{result}\label{res: d_graph_results}
Given $\epsilon>0$ and $n\in\integers_+$, let $n=2\tn$. Among all the gossiping algorithms $\mathcal{A}(P)$ with a symmetric $P$ on the barbell graph, $K_\tn-K_\tn$, 
randomized ER gossiping leads to $T_{ave}(\varepsilon, P^r)=\Theta(n^2\log(1/\varepsilon))$, which is optimal with respect to $\varepsilon$ and $n$, and cannot be improved.

In a more general setting, let $\cG$ be a connected graph with diameter $\mathcal{D}\in\integers_+$. The $\varepsilon$-averaging time of $\mathcal{A}(P^r)$ satisfies
$$
T_{ave}(\varepsilon, P^r) = \mathcal{O}(
\mathcal{D}n^3)\log(\epsilon^{-1}).
$$
\end{result}}
\begin{remark} The $\varepsilon$-averaging time of randomized gossiping with lazy Metropolis weights\footnote{For lazy Metropolis weights see \eqref{Metr-weights} and the paragraph after.} on any graph is $\cO(n^3 \log(1/\varepsilon))$; while, for the barbell graph, Metropolis weights perform similar to uniform weights; both require $\Theta(n^3 \log(1/\varepsilon))$ time which can be improved to $\Theta(n^2 \log(1/\varepsilon))$ by 
ER gossiping.
\end{remark}

\begin{remark}
If the diameter $\mathcal{D}\leq 11$, our bounds for 
ER gossiping improve upon that of the randomized gossiping with lazy Metropolis weights 
by a (small) constant factor (see Remark \ref{remark-metropolis-compare}).
Note $\mathcal{D}=3$ for barbell graphs and $\cD\leq 11$ is also reasonable for mid-size \textit{small-world} graphs which are random graphs that arise frequently in real-world applications \cite{ContSmallWorld}. 
{For instance, Cont \emph{et al.} \cite{ContSmallWorld} show that 
the diameter $\mathcal{D}$ of the randomized community-based small-world graphs admits $2\log(n)$ upper bound almost surely; hence, for these graphs 
$\mathcal{D}\leq 11$ almost surely for $n\leq 240$.
Indeed, we empirically observe that randomly generated 
small-world graphs 
with parameters $n=\{5k:k=1,\ldots,5\}$ and $m=\lfloor 0.2(n^2-n)\rfloor$ using the methodology described in the numerical experiments in Section~\ref{sec:numerics-consensus}} satisfy $\mathcal{D} \leq 5$ on average over $10^4$ independent and identically distributed (i.i.d.) samples.

\end{remark}

\vspace*{-5mm}
\section{Proofs of Main Results}
\label{sec-theoretical guarantees}


In order for 
both uniform and ER gossiping methods to have the same expected number of node wake-ups in a given time period, one should have $r_i=r=2(n-1)/n$ for $i\in\cN$ within 
the uniform gossiping model --recall that $r_i=\sum_{j\in\cN_i}R_{ij}$ for $i\in\cN$ for ER gossiping; hence, the rate of both Poisson processes will be the same, i.e., $\sum_{i\in\cN}r_i=2(n-1)$. We note that the number of clock ticks $k\in\integers_+$ can be converted to \emph{absolute time} easily with standard arguments (simply dividing $k$ by $\sum_{i\in\cN}r_i$ to get the expected time of the $k$-th tick), e.g., see \cite[Lemma 1]{boyd2006randomized}. This allows us to use the number of iterations (clock ticks) to compare asynchronous algorithms.

{It can be easily 
verified that 
for a given $\cA(P)$, the expected iteration matrix defined in \eqref{def-Wp} satisfies
\begin{align}
\label{eq:W-P_identity}
    \overline{W}_P=I-\tfrac{1}{2}D+\tfrac{1}{2}(P+P^\top),
\end{align}
where $D$ is a diagonal matrix with $i$-th entry $D_i\triangleq \sum_{j\in\cN_i}(P_{ij}+P_{ji})$. Note $W_{ij}$ defined in Section~\ref{sec:gossiping} is a doubly stochastic, non-negative and weakly diagonally dominant matrix for all $i\in\cN$ and $j\in\cN_i$; therefore, $\overline{W}_P$, which is a convex combination of $W_{ij}$ matrices, is also a doubly stochastic, non-negative and weakly diagonally dominant matrix. It follows then from the Gershgorin's Disc Theorem (see e.g. \cite{golub2012matrix}) that all the eigenvalues of $\overline{W}_P$ are non-negative. Moreover, since $\overline{W}_P$ is a non-negative doubly stochastic matrix, its largest eigenvalue $\lambda_n(\overline{W}_P) = 1$.}
 Plugging in $P^u$ and $P^r$ for $P$ in this identity respectively leads immediately to the following result. 
\begin{lemma}\label{lemma:WP-formula} The matrices $\overline{W}_{P^r}=\mathbb{E}_{P^r} [W_{ij}] $ and $\overline{W}_{P^u}=\mathbb{E}_{P^u} [W_{ij}]$ satisfy the following identities:
{\small
\begin{eqnarray*}
    \overline{W}_{P^u}=I-\frac{1}{2}D^u+\frac{P^u+(P^{u})^\top}{2}, \quad
    \overline{W}_{P^r}=I-\frac{1}{2}D^r+P^r, 
\end{eqnarray*}}%
where $D^u$ and $D^r$ are diagonal matrices satisfying $[D^u]_{ii} \triangleq \sum_{j\in\cN_i}(P^u_{ij}+P^u_{ji})$, $[D^r]_{ii}=\frac{1}{(n-1)}R_{i}$ where $R_i \triangleq \sum_{j\in\cN_i}R_{ij}$. 
\end{lemma}

Recall the definition of $T_{ave}(\varepsilon,P)$ given in \eqref{eq:T_ave}, i.e.,  $\varepsilon$-averaging time of an asynchronous gossiping algorithm $\mathcal{A}(P)$ characterized by 
a probability matrix $P$. 
According to Theorem~\ref{thm-gossip-compare}, 
to compare 
uniform and ER gossiping methods introduced in Section~\ref{sec:prelim}, it is sufficient to estimate the second largest eigenvalues of $\overline{W}_{P^r}$ and $\overline{W}_{P^u}$ and compare them. In the rest of this section, we discuss estimating the second largest eigenvalues of $\overline{W}_{P^r}$ and $\overline{W}_{P^u}$ based on the notions of graph conductance and hitting times when the eigenvalues are not readily available in closed form. We will also discuss some examples 
for which we can explicitly compute the eigenvalues.

It is worth emphasizing that since the matrices $\overline{W}_{P^r}$ and $\overline{W}_{P^u}$ are symmetric and doubly stochastic, 
they can both be viewed as the probability transition matrix of a \emph{reversible Markov Chain} on the graph $\cG$, both with a uniform stationary distribution.  
We saw that depending on the type of randomized gossiping, the sampling probabilities of the bottleneck edge can differ significantly --by a factor of $\Theta(n)$ on barbell graphs implied by \eqref{eq-bottleneck-uniform-proba} and \eqref{eq-bottleneck-resist-proba}. A similar effect can also be observed for the Markov chains defined by the transition probability matrices $\overline{W}_{P^u}$ and $\overline{W}_{P^r}$. In fact, by an explicit computation based on Lemma~\ref{lemma:WP-formula} (see Lemma \ref{lemma: Values} for details), we get
    $$[\overline{W}_{P^u}]_{i^*j^*}= 
    \frac{2}{n^2}, \quad [\overline{W}_{P^r}]_{i^*j^*}
    =\frac{1}{2(n-1)}.$$
That is, the probability of moving from one complete subgraph to the other is significantly larger (by a factor of $\Theta(n)$) for the Markov chain corresponding {to} 
$\overline{W}_{P^r}$ than that of the chain with $\overline{W}_{P^u}$. Intuitively speaking, this fact allows the 
ER-based chain to traverse between the 
complete subgraphs faster when $n$ is large, 
leading to faster averaging over the nodes. This will be formalized and proven in the next subsection, where we study gossiping algorithms over \emph{barbell} and \emph{$c$-barbell graphs}. 
\subsection{Proof of Result\;\ref{res: c-barbell graph} via conductance-based analysis}
\label{Subsec: Conductance}
Probability transition matrices on graphs have been studied well; 
in particular, there are some combinatorial techniques to bound their eigenvalues based on \emph{graph conductance}~\cite{aldous-fill-2014} as well as some algebraic techniques that allow one to compute all the eigenvalues explicitly exploiting symmetry groups of a graph \cite{BarbellBoyd} as we shall discuss in Section~\ref{sec:spectral}. 

The notion of graph conductance is tied to a transition matrix $W$ over a graph which corresponds to a reversible Markov chain admitting an arbitrary stationary 
distribution $\pi$. 
It can be viewed as a measure of how hard it is for the Markov chain to go from a subgraph to its complement in the worst case. 

The notion of graph conductance allows us to provide bounds on the mixing time of the corresponding Markov chain as we discuss below.
\begin{defi}[Conductance]
 Let 
 $W$ be the transition matrix of a reversible Markov chain\footnote{That is $\pi_iW_{ij}=\pi_jW_{ji}$ for all $i,j \in \mathcal{N}$.} on the graph $\cG$ with a stationary distribution $\pi=\{\pi_i\}_{i=1}^n$. The conductance $\Phi$ is defined as

 \BEQ\label{eq: conductance}
\Phi(W)=\min_{\substack{S\subset \mathcal{N}:\sa{S,S^c\neq\emptyset}}} \frac{\sum_{i\in S,j\in S^c}\pi_i W_{ij}}{\min\{\pi(S),\pi(S^c)\}}
 \EEQ
  
 where $\pi(S)\triangleq\sum_{i\in S}\pi_i$. 
 \end{defi}

Given a transition matrix $W$, the relation between conductance $\Phi(W)$ and the second largest eigenvalue $\lambda_{n-1}(W)$ is well-known and given by the \emph{Cheeger inequalities}:
\beq \label{ineq-Cheeger}
1- 2\Phi(W) \leq \lambda_{n-1}(W) \leq 1-\Phi^2(W),
\eeq
--see, e.g., \cite[Proposition 6]{DiaconisBounds}. Therefore, larger conductance leads to faster averaging, i.e., shorter $T_{ave}(\varepsilon,P)$, in light of Theorem~\ref{thm-gossip-compare}. 
In particular, we can get lower and upper bounds on the averaging time for both 
uniform and ER gossiping methods using the Cheeger's inequality. We study the performance bounds for these gossiping algorithms over \emph{c-barbell graphs}; 
and our next result shows $\Theta(n)$ improvement on the conductance of effective resistance-based transition probabilities $\overline{W}_{P^r}$ compared to uniform probabilities $\overline{W}_{P^u}$ on a c-barbell graph with $n=c\tn$ nodes. 

\begin{proposition} \label{prop: c-barbell graph} Given $\tilde{n},c\in\integers_+$ such that $c\geq 2$, consider the two Markov chains on the $c$-barbell graph with $n=\tn c$ nodes defined by the transition matrices $\overline{W}_{P^u}$ and $\overline{W}_{P^r}$. 
Let $c_*=\big(\lfloor \frac{c}{2}\rfloor\big)^{-1}$. The 
conductance values are given by
  \begin{align}\label{Cond: c-barbell}
    \Phi(\bar{W}_{P^u})=\frac{c_*}{c\tn^3},\;\;\;\; \Phi(\bar{W}_{P^r})=\frac{c_*}{2\tn(c\tn-1)}.
  \end{align} 

\end{proposition}
\begin{remark}
\label{rem:barbel-conductance}
Since a barbell graph $K_\tn - K_\tn$ is a special case of a c-barbell graph with $c=2$ and $n=2\tn$, Proposition~\ref{prop: c-barbell graph} implies that $\Phi(\overline{W}_{P^u})=\frac{4}{n^3}$ and  
    $\Phi(\overline{W}_{P^r})=\frac{1}{n(n-1)}$.
\end{remark}
Given the transition matrix $W$, by taking the logarithm of the Cheeger inequalities in~\eqref{ineq-Cheeger}, for $\Phi(W)\leq 1/2$, we obtain
{\small
\beq
\label{eq:lambda_bound}
    -\log(1-\Phi^2(W)) \leq \log (\lambda_{n-1}^{-1}(W)) \leq -\log(1-2\Phi(W)). 
\eeq}%
Then, choosing $W=\overline{W}_{P^u}$ and $W=\overline{W}_{P^r}$ above,  applying Theorem \ref{thm-gossip-compare} and Proposition~\ref{prop: c-barbell graph}
and noting $-\log(1-x) \approx x$ for $x$ close to 0, leads to the 
lower and upper bounds on the averaging time of 
uniform and ER gossiping algorithms as shown in Result~\ref{res: c-barbell graph} {
of our main results section} (Section \ref{sec: Main Results}). \mg{In the Supplementary Material, we also studied the tightness of our conductance bounds \eqref{eq:lambda_bound} numerically on the $c$-barbell graphs to show that our bounds are reasonable. In particular, we observe that our lower bounds gets tighter as the number of nodes, $n$, increases \bc{on c-barbell graphs}.}
  

Although this analysis is also applicable to other graphs with low conductance, it does not typically lead to tight estimates, i.e., the lower and upper bounds do not match 
in terms of their dependency on $n$. In the next section, we show that for the case of barbell graphs, we get tight estimates on the averaging time by computing the eigenvalues of the averaging matrices $\overline{W}_{P^r}$ and $\overline{W}_{P^u}$ \emph{explicitly}. More precisely, we will show in Proposition \ref{prop: eig-order-imprv} that the lower bounds in 
\eqref{bound: c_barbell_uni}--\eqref{bound:c_barbell_res} are tight for $c=2$ in the sense that $T_{ave}(\varepsilon,P^u ) = \Theta(n^3)$ and $T_{ave}(\varepsilon,P^r ) = \Theta(n^2)$ and the effective resistance-based averaging is faster by a factor of $\Theta(n)$ which will imply Result \ref{res: eig-order-imprv}.

\subsection{Proof of Result\;\ref{res: eig-order-imprv} via spectral analysis}
\label{sec:spectral}
\begin{figure}[h!]
    \centering
    \includegraphics[width=0.60 \linewidth]{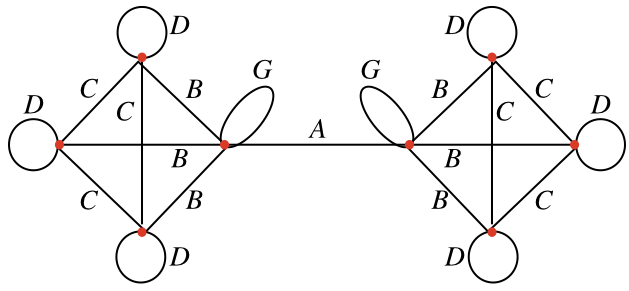}
    \caption{\label{fig-edge-weighted-barbell}An edge-weighted barbell graph $K_{\tilde{n}}-K_{\tilde{n}}$ with edge weights $A,B,C,D,G > 0$ for $\tilde{n}=4$.}
\end{figure}

Eigenvalues of probability transition matrices defined on  barbell graphs are studied in the literature. Consider the \emph{edge-weighted} barbell graph $K_\tn-K_\tn$ with $n=2\tn$ nodes, where $w=[w_{ij}]_{(i,j)\in\cE}$ is the vector of edge weights that have positive entries. Suppose each node has a self-loop, e.g., see Fig.~\ref{fig-edge-weighted-barbell}. Let $(i^*,j^*)$ be the edge that connects the two complete subgraphs. The result \cite[Prop. 5.1]{BarbellBoyd} gives an explicit formula for the eigenvalues of a probability transition matrix $W$ with transition probabilities proportional to edge weights, i.e., $W_{ij} = w_{ij} / \sum_{j\in\cN_i} w_{ij}$ where $w_{ij}$ satisfy the following assumptions: $w_{i^*i^*}=w_{j^*j^*}=0$, $w_{i^*j^*}=A$, $w_{i^*j}=w_{j^*i}=B$ for all $j\in\cN_{i^*}\setminus\{j^*\}$ and $i\in\cN_{j^*}\setminus\{i^*\}$, $w_{ij}=C$ for all $(i,j)$ in each $K_{\tilde{n}}$ such that $i\neq j$ and $i,j \notin \{i^*,j^*\}$, and $w_{ii}=D$ for $i\in\cN\setminus\{i^*,j^*\}$ for some $A,B,C,D>0$. 
Note we cannot immediately use this result to compute the eigenvalues of the transition matrices $\overline{W}_{P^r}$ and $\overline{W}_{P^u}$ defined in Lemma~\ref{lemma:WP-formula}. Mainly because
all the diagonal entries 
of $\overline{W}_{P^r}$ and $\overline{W}_{P^u}$ being strictly positive breaks the $w_{i*i*}=w_{j*j*}=0$ assumption of \cite[Prop. 5.1]{BarbellBoyd}.
In Proposition~\ref{prop: gen. eigenvalue}, we adapt 
\cite[Prop. 5.1]{BarbellBoyd} to our setting with some minor modifications to allow $w_{i^*i^*}=w_{j^*j^*}=G$ for any $G>0$ so that it becomes applicable to $\overline{W}_{P^r}$ and $\overline{W}_{P^u}$. The proof of Proposition~\ref{prop: gen. eigenvalue}, {provided in the Supplementary Material}, is similar to the proof of \cite[Prop. 5.1]{BarbellBoyd} and is based on exploiting the symmetry properties of the weighted barbell graph as described above --illustrated in Figure \ref{fig-edge-weighted-barbell} for $\tilde{n} = 4$.
\begin{proposition}[Generalization of Proposition~5.1 in~\cite{BarbellBoyd}] \label{prop: gen. eigenvalue} Consider the edge-weighted barbell graph $K_\tn-K_\tn$ with $n=2\tn$ nodes. Let $(i^*,j^*)$ be the edge that connects the two complete subgraphs. Assume that weights are of the form $w_{i^*i^*}=w_{j^*j^*}=G$, $w_{i^*j^*}=A$, $w_{i^*j}=w_{j^*i}=B$ for all $j\in\cN_{i^*}\setminus\{j^*\}$ and $i\in\cN_{j^*}\setminus\{i^*\}$, $w_{ij}=C$ for all $(i,j)$ in each $K_{\tilde{n}}$ such that $i\neq j$ and $i,j \notin \{i^*,j^*\}$, and $w_{ii}=D$ for $i\in\cN\setminus\{i^*,j^*\}$ for some $A,B,C,D,G>0$. Consider the transition matrix $W$ associated to this graph with entries $W_{ij} = w_{ij} / \sum_{j\in\cN_i} w_{ij}$, then the eigenvalues of $W$ are 
\begin{itemize}
    \item $\lambda_{a}\triangleq 1$ with multiplicity one,
    \item $\lambda_{b} \triangleq -1+\frac{A+G}{A+G+E}+\frac{F}{F+B}$ with multiplicity one,
    \item $\lambda_{c} \triangleq \frac{D-C}{B+F}$ with multiplicity $n-4$,
    \item $\lambda_{\pm} \triangleq \frac{1}{2}\Big(\frac{F}{B+F}+\frac{G-A}{A+E+G}\,\pm\,\sqrt{S} \Big)$,
\end{itemize}
where $E \triangleq (\tn-1)B$, $F\triangleq D+(\tn-2)C$ and $S \triangleq\big(\frac{F}{B+F}+\frac{G-A}{A+E+G}\big)^2-\frac{4(FG-BE-AF)}{(B+F)(A+E+G)}$. 
\end{proposition}
Based on this 
result, in Proposition~\ref{prop: eig-order-imprv} we characterize the second largest eigenvalue of the transition matrices $\overline{W}_{P^u}$ and $\overline{W}_{P^u}$ 
-- the proof can be found in the appendix.
\begin{proposition}\label{prop: eig-order-imprv} 
Consider Markov chains on the barbell graph $K_{\tn}-K_{\tn}$ with transition matrices $\overline{W}_{P^r}$ and $\overline{W}_{P^u}$. The second largest eigenvalues of these matrices are given by
\begin{eqnarray*}
    \lambda_{n-1}(\overline{W}_{P^r}) = 1 - \Theta(\frac{1}{n^2}),\;\;\;\; 
   \lambda_{n-1}(\overline{W}_{P^u}) = 1 - \Theta(\frac{1}{n^3}).
\end{eqnarray*}
\end{proposition}
Result~\ref{res: eig-order-imprv} follows as a direct consequence of Proposition \ref{prop: eig-order-imprv} and Theorem~\ref{thm-gossip-compare}. Thus, we establish that that averaging time with resistance weights is $\Theta(n)$ faster on a barbell graph.

\subsection{Proof of Result\;\ref{res: d_graph_results} via hitting and mixing times}
Before giving a formal definition of the $\varepsilon$-mixing time, we introduce the total variation (TV) distance between two probability measures $p$ and $q$ defined on the set of nodes $\cN = \{1,2,\dots,n\}$. 
TV distance between $p$ and $q$ is defined as $\| p - q \|_{TV} \triangleq  \| p - q\|_1/2.$
\noindent Given a Markov chain $\mathcal{M}$ with a probability transition matrix $W$ and stationary distribution $\pi$, $\varepsilon$-mixing time is a measure of how many iterations are needed for the probability distribution of the chain to be $\varepsilon$-close to the stationary distribution 
in the TV distance. {A related notion 
is the hitting time which is a measure of how fast the Markov chain travels between any two nodes.} 
\begin{defi} {(Mixing time and hitting times)}
Given $\varepsilon>0$ and a Markov chain with probability transition matrix $W$ and stationary distribution $\pi$, the $\varepsilon$-mixing time is defined as 
    $$T_{mix}(\varepsilon, W) \triangleq \inf_{k\geq 0}  \left\{ \sup_{p\geq 0: \|p\|_1 = 1} \| (W^k)^\top p - {\pi} \|_{TV} \leq \varepsilon \right\},$$
{and the hitting time $\htw$ is the  expected number of steps until the Markov chain reaches $j$ starting from $i$.}
\end{defi}
Mixing-times and averaging times are closely related. In fact, given probability transition matrix $W$, it is known that $T_{ave}(\varepsilon,W)$ and $T_{mix}(\varepsilon, \tilde{W})$ admit the same bounds up to $n \log n$ factors \cite[Theorem 7]{boyd2006randomized} for $\tilde{W} = \frac{I + W}{2}$.\footnote{Note \cite[Theorem 7]{boyd2006randomized} uses absolute time whereas we used number of node wake-ups to define $\epsilon$-averaging and $\epsilon$-mixing times; therefore, we multiplied $\log(n)$ factor in~\cite[Theorem 7]{boyd2006randomized} by $\sum_{i\in\cN}r_i=2(n-1)$ to convert absolute times to number of node wake-ups.} Hence, designing algorithms with a smaller mixing time, often leads to better algorithms for distributed averaging (see also \cite{shah-book}). It is also known that mixing time is closely related to hitting times \cite[Theorem 1.1]{peresmixing}. 
 
 Next, we show the first part of Result~\ref{res: d_graph_results}, i.e.,  $T_{ave}(\varepsilon,P^r)=\Theta(n^2 \log(1/\varepsilon))$ is optimal among all $\cA(P)$ with a symmetric $P$. Note 
 $P$ is symmetric implies that it is doubly stochastic. For large $n$ and doubly stochastic $P$, by \cite[Corollary 1]{boyd2006randomized}, we have $T_{ave}(\varepsilon,P) = \Theta\left( \frac{n \log(1/\varepsilon)}{1-\lambda_{n-1}(P)}\right)$. \bc{On the other hand, Roch proved in \cite[Section 3.3.1]{Bound_Roch} that any \sa{symmetric} doubly stochastic $P$ matrix on the barbell graph with $n$ nodes satisfies the bound $\frac{1}{1-\lambda_{n-1}
(P)}=\Omega(n)$. Inserting this estimate into the expression for the averaging time, we obtain $T_{ave}(\varepsilon,P) = \Omega\left(n^2 \log(1/\varepsilon)\right)$ for any $\mathcal{A}(P)$ with symmetric $P$ on barbell graphs}.

\mg{We conclude that the averaging time of the ER-based gossiping on the barbell graph, which satisfies $T_{ave}(\varepsilon,P^r ) = \Theta(n^2\log(1/\varepsilon))$ by Proposition \ref{prop: eig-order-imprv} and Theorem \ref{thm-gossip-compare}, is optimal with respect to its dependency to $n$ and $\varepsilon$ among all symmetric choices of the $P$ matrix.}

Next, given any connected graph $\mathcal{G}$, we obtain a bound on the second largest eigenvalue of the $\overline{W}_{P^r}$ and show that the averaging time with effective
resistance weights  $T_{ave}(\varepsilon,P^r) = \mathcal{O}\left(\mathcal{D} n^3 \log(1/\varepsilon)\right)$ where 
$\mathcal{D}$ is the diameter of the graph.
\begin{theorem}\label{thm: d_graph_results} Let $\cG$ be a graph with diameter $\mathcal{D}$. The second largest eigenvalue of 
$\overline{W}_{P^r}$ 
satisfies
    $\lambda_{n-1}(\overline{W}_{P^r}) \leq 1 - \frac{1}{6\mathcal{D} n^3 }.$
\end{theorem}
 \begin{proof} 
 It follows from our discussion in Section~\ref{sec-theoretical guarantees} that $\overline{W}_{P^r}$ is non-negative and doubly stochastic (see the paragraph before Lemma \ref{lemma:WP-formula}). 
Therefore, for analysis purposes, we can interpret $\overline{W}_{P^r}$ as the transition matrix of a Markov chain $\mathcal{M}$ whose stationary distribution $\pi$ is the uniform distribution. Our analysis is based on relating the eigenvalues of $\overline{W}_{P^r}$ matrix to the hitting times of the Markov chain $\mathcal{M}$ where we follow the proof technique of \cite[Lemma 2.1]{olshevsky2015linear}. 
By Lemma \ref{lem-hitting-ub} from the appendix, we get
	$ \hwprij \leq n \frac{2(n-1)}{R_{ij}}$ if $j\in \cN_i.$
For any graph, it is also known that\footnote{This follows directly from the Rayleigh's monotonicity rule \cite{doyle1984random}  which says that if an edge is removed from a graph, effective resistance on any edge can only increase. Therefore, the complete graph provides a lower bound for $R_{ij}$ where $R_{ij}=2/n$ (see also {\cite{chandra1996electrical}}).}	
    $\min_{i,j} R_{ij} \geq \frac{2}{n}.$
Therefore, for any neighbors $i$ and $j$, 
$ \hwprij \leq n^2 (n-1).$
For any two vertices $i$ and $j$ not necessarily neighbors, $i\neq j$, let $v_0(=i), v_1, \dots, v_\ell(=j)$ be the shortest path connecting $i$ and $j$. Then, by the subadditivity property of hitting times, for any $i,j\in\cN$, we obtain
    $\hwprij \leq \ell n^2 (n-1) 
    \leq \mathcal{D} n^2 (n-1)$.
It follows from an analysis similar to \cite{levin2009markov} that 
\begin{eqnarray} 
\begin{footnotesize}
T_{mix}(\frac{1}{8}, \overline{W}_{P^r}) \leq 8 \max_{i,j\in \{1,\dots,n\}}\hwprij  + 1  
 \leq 8\mathcal{D} n^3. \label{ineq-mix-time-ub}
\end{footnotesize}
\end{eqnarray}
From \cite[eqn. (12.12)]{levin2009markov}, we also have
 $$ T_{mix}(\frac{1}{8}, \overline{W}_{P^r}) \geq \left(\frac{1}{1-\lambda_{n-1}(\overline{W}_{P^r})} -1 \right) \log(4).$$
Combining this with the estimate \eqref{ineq-mix-time-ub} implies directly $\lambda_{n-1}(\overline{W}_{P^r}) \leq 1 - \frac{1}{6\mathcal{D} n^3},$
which proves the claim.
\end{proof}
{\bf Metropolis vs ER gossiping:}
Given a connected $\cG=(\cN,\cE)$, suppose there are no self-loops, i.e., $(i,i)\not\in\cE$ for $i\in\cN$. Uniform weights $p^{u}_{j|i}=\frac{1}{d_i}$ can result in slow mixing on some graphs such as the barbell graph (see  Proposition \ref{prop: eig-order-imprv}) or other graphs like lollipop graphs \cite{aldous-fill-2014} which have both high degree and low degree nodes together. A popular alternative to uniform weights $\{p^{u}_{j|i}\}_{j\in\cN_i}$ for $i\in\cN$ is the Metropolis weights defined as 
\begin{equation} \label{Metr-weights}
M_{ij} \triangleq \begin{cases} 
    \frac{1}{\max(d_i, d_j)} & \mbox{if } (i,j) \in \mathcal{E},\\
    1-\sum_{j\in \bc{\mathcal{N}_i\setminus i}} \frac{1}{\max(d_i, d_j)} & \mbox{if } i=j, \\
    0                        & \mbox{else}.
\end{cases} 
\end{equation} 

\mg{Let $M$ denote the matrix whose entries are the Metropolis weights $M_{ij}$.}  
\mg{The weights determined by the matrix ${\widetilde{M}} \triangleq \frac{I + M}{2}$ are also
popular in the distributed optimization practice \cite{olshevsky2015linear} which is referred to as the lazy version \mg{of the Metropolis weights}}. The matrix {$\widetilde{M}$} is symmetric and positive semi-definite, unlike the matrix $M$ which may have negative eigenvalues that can be close to $-1$ (therefore, it can be problematic for the convergence of distributed algorithms, see e.g. \cite{shi2015extra}). Combined with uniform wake-up of nodes, this leads to the following wake-up probabilities for the Metropolis weights based system:
${P_{ij}^{\bc{\widetilde{m}}}} \triangleq \frac{1}{n}\widetilde{M}_{ij},$
and the associated matrix
$\overline{W}_{P^{\bc{\widetilde{m}}}} \triangleq  \mathbb{E}_{P^{\bc{\widetilde{m}}}} [\bc{W^{(i,j)}}] = \sum_{ij} P_{ij}^{\widetilde{m}} \bc{W^{(i,j)}}. $
In particular, for any connected graph $\cG=(\mathcal{N},\mathcal{E})$ with $n$ nodes, we have the following guarantees from \cite[Lemma 2.1]{olshevsky2015linear} on the lazy Metropolis weights: 
  \begin{equation}\label{ineq-metropolis-eig}\max_{i,j \in \{1,2,\dots,n \} } \htmij  \leq 12n^2, \quad \lambda_{n-1} (\widetilde{M}) \leq 1- \frac{1}{71 n^2}.
  \end{equation}
By \eqref{eq:W-P_identity}, we have also 
$   \overline{W}_{P^{\widetilde{m}}}=(1-\frac{1}{n}) I + \frac{1}{n}\widetilde{M}.
$
Therefore, from \eqref{ineq-metropolis-eig}, we get the bound 
$
 \lambda_{n-1}(\overline{W}_{P^{{\widetilde{m}}}}) \leq 
 1 - \frac{1}{71n^3}
 $,
 for any connected graph $\cG$. Therefore, we conclude from Theorem \ref{thm-gossip-compare} that the $\varepsilon$-averaging time of Metropolis weights-based gossiping on any graph is $\cO(n^3 \log(1/\varepsilon))$ -- again using the fact that $-\log(1-x)\approx x $ for $x$ close to $0$.
 That said, for barbell graphs, Metropolis weights perform similar to uniform weights; both require $\Theta(n^3 \log(1/\varepsilon))$ time which is improved by the effective resistance-based weights to $\Theta(n^2 \log(1/\varepsilon))$. This completes the proof of Result \ref{res: d_graph_results}.
 
 \begin{remark}\label{remark-metropolis-compare} Comparing the inequalities $\lambda_{n-1}(\overline{W}_{P^r}) \leq 1 - \frac{1}{6\mathcal{D} n^3}$ and $\lambda_{n-1}(\overline{W}_{P^{\bc{\widetilde{m}}}}) \leq 
 1 - \frac{1}{71n^3}$, we see that for $\mathcal{D} \leq 11$, the upper bound 
 on $\lambda_{n-1}(\overline{W}_{P^r})$ will be smaller than the upper bound for $\lambda_{n-1}(\overline{W}_{P^{\bc{\widetilde{m}}}})$. Therefore, performance bounds obtained on the $\varepsilon$-averaging time through Theorem \ref{thm-gossip-compare} for ER weights will be better than those of Metropolis weights by a (small) constant factor for $\mathcal{D} \leq 11$.
 \end{remark}
\vspace{-4mm}
\section{Numerical Experiments}\label{sec-numerical}
In this section, we demonstrate the benefits of using effective resistances for solving the consensus problem and also within DPGA-W~ \cite{aybat2018distributed} and EXTRA~\cite{shi2015extra} algorithms for consensus optimization.

\subsection{Consensus exploiting effective resistances}
\label{sec:numerics-consensus}
Gossiping algorithms have been studied extensively 
and there have been a number of approaches \cite{shah-book,GossipTsitsiklis,GossipZanaj,GossipKermarrec,GossipKempe,GossipFagnani,GossipEstrin}.
In light of Theorem \ref{thm-gossip-compare}, among all the algorithms $\mathcal{A}(P)$ with a symmetric $P$, the matrix $P^{opt}$ that minimizes the second largest eigenvalue, i.e. $\lambda_{n-1}(\overline{W}_{P})$, is the fastest. The gossiping algorithm $\mathcal{A}(P^{\text{opt}})$ with optimal choice of the probability matrix $P^{\text{opt}}$ is called the Fastest Mixing Markov Chain (FMMC) in the literature \cite{Boyd03fastestmixing}. 
In \cite{boyd2006randomized}, Boyd et al. propose a distributed subgradient method to compute the matrix $P^{\text{opt}}$.
This method requires a decaying step size and computation of the subgradient of the objective $\lambda_{n-1}(\overline{W}_P)$ with respect to the decision variable $P$ at every iteration which itself requires solving a consensus problem at every iteration. This can be expensive in practice in terms of average number of communications required per node, and its convergence to $P^{opt}$ can be slow with at most sublinear convergence rate \cite{boyd2006randomized}.
{In contrast, 
ER probabilities $P^r$ are optimal for some graphs (such as the barbell graph, see Result \ref{res: d_graph_results}) and can be computed efficiently with the normalized D-RK algorithm (see the Supplementary Material) which admits linear convergence guarantees. Therefore, ER weights can serve as a computationally efficient alternative to optimal weights for consensus. For illustrating this point, 
we compare communication requirements per node for ER gossiping and FMMC on barbell and small-world graphs. This comparison consists of two stages: $(i)$ pre-computation stage (where the probability matrices $P^r$ and $P^{\text{opt}}$ are computed up to a given tolerance) $(ii)$ asynchronous consensus stage (where we run ER and FMMC with probability matrices $P^r$ an $P^{\text{opt}}$ obtained from the previous stage to solve a consensus problem). 
}

First, we implement subgradient method with decaying step size $\alpha_k = R/k$ from \cite{boyd2006randomized} where $R$ is tuned to the graph to achieve the best performance and stop the computation of matrix of FMMC at step $k$ if the iterate $P_k$ satisfies $\frac{||P_k-P^{opt}||_F}{||P^{opt}||_F} \leq \epsilon_1$ where $\epsilon_1$ is the given precision level.\footnote{The optimal probability matrix $P^{opt}$ which serves as a baseline in the stopping criterion is estimated accurately by solving the semi-definite program (SDP) \cite[eqn. (53)]{boyd2006randomized} directly using the CVX software \cite{cvx} with a centralized method and computations required to solve this SDP is not counted as a part of the communication cost we report for FMMC in Tables \ref{Table: Barbell}-\ref{tab:Comparison-SmallWorld}.} 
Similarly, we compute $\mathcal{L}^{\dag}$ for ER and stop the normalized D-RK algorithm when the iterate $X^k$ at step $k$ satisfies $\frac{\|X^k-\cL^\dag\|_F}{\|
{\cL}^\dag\|_F} \leq \epsilon_1$. Since the distributed subgradient method of \cite{boyd2006randomized} is based on synchronous computations, we also implemented the normalized D-RK algorithm with synchronous computations for fairness of comparison. {We define the \emph{communication}  for a node as a contact with its neighbour either to compute an average of their state vectors or to update the matrix $P_k$ at 
any iteration.}

We compared both of the algorithms based on their communication performances on stage-i an stage-ii. In particular, we considered the number of communications required per node to obtain the matrix $P_k$ for ER and FMMC at stage-i and at stage-ii, we generated 1000 instances of $y_i^0$ to start consensus and compare the average number of communications per node required to achieve $y_i^k$ satisfying $\frac{||y^k-\bar{y}||}{||\bar{y}||}\leq\epsilon_2$ where $\epsilon_2$ is the tolerance level.

{For the barbell graph, the initial state vector $y_i^0$ for consensus is sampled from the normal distribution $\textbf{N}(500,10)$ if $i\in \mathcal{N}_{L}$ and from $\textbf{N}(-500,10)$ if $i\in \mathcal{N}_R$ where tolerance levels are set to be $\epsilon_1=\epsilon_2=0.01$. We also compare ER and FMMC on small-world graphs while the number of nodes $n$ is varied with an edge density $\frac{2m}{n^2-n}\approx 0.4$ {where $m$ is the total number of edges}. On small-world graphs we generated $1000$ instances of $y_i^0$ drawn from $\textbf{N}(0,100)$ and stopped algorithms whenever tolerance levels $\epsilon_1=\epsilon_2=0.05$ are obtained or the number of communications per node exceeded $10^6$. }

\begin{table}[!ht]
\centering
\begin{tabular}{|r|l|l|l|}
\hline
\multicolumn{1}{|l|}{Graph} & Method & \begin{tabular}[c]{@{}l@{}l@{}l@{}} Comm.\\per node\\(stage-$i$)\\\end{tabular} & \begin{tabular}[c]{@{}l@{}l@{}}Comm. \\ per node \\ \big(\begin{tabular}[c]{@{}l@{}} stage-$ii$\end{tabular}\big) \end{tabular} 
\\ \hline
$K_5-K_5$ & ER & 2.9  $\times 10^3$  & 81 \\ 
 & FMMC & 1.28  $\times 10^5$ & 65 \\ \hline
$K_{10}-K_{10}$ & ER & 8.4 $\times 10^4$ & 198 \\ 
 & FMMC & 3.93  $\times 10^5$ & 130 \\ \hline
$K_{20}-K_{20}$ & ER & 2.6 $\times 10^6$ & 433  \\ 
 & FMMC & 6.4 $\times10^6$ & 251 \\ \hline
$K_{25}-K_{25}$ & ER & 7.9  $\times 10^6$ & 566 \\ 
 & FMMC & $>10^{7}$ & 287 \\ \hline
\end{tabular}
\caption{
FMMC vs ER on the barbell graph.
}
\label{Table: Barbell}
\end{table}

Results for both of the graphs are reported in Tables \ref{Table: Barbell} and \ref{tab:Comparison-SmallWorld} in which we compare the average communication per node in the pre-computation (stage-$i$) and in the consensus computation (stage-$ii$) where results are averaged over 1000 runs.
On barbell graph, we observe that FMMC requires less communications at the second (consensus) stage as expected (as FMMC is based on the optimal matrix $P^{opt}$), but in terms of total communications (stage-$i$ + stage-$ii$) ER outperforms FMMC. In the case of small-world graphs, computation of $P^{opt}$ 
exceeded the maximum communication limit which caused FMMC to perform worse than ER in stage-$ii$ (since the stage-$i$ solution is not a precise approximation of $P^{opt}$ anymore). We can say that ER performs better than FMMC in terms of total communications for both graph types.

\begin{table}[ht!]
\centering
\begin{tabular}{|r|l|l|l|}
\hline
\multicolumn{1}{|l|}{Graph} & Method & \begin{tabular}[c]{@{}l@{}l@{}}Comm.\\per node\\(stage-$i$)\end{tabular} & \begin{tabular}[c]{@{}l@{}l@{}}Comm. \\ per node \\ \big(\begin{tabular}[c]{@{}l@{}}stage-$ii$\end{tabular}\big) \end{tabular} 
\\ \hline
$n=5$ & ER & 6.4 & 41 \\
 & FMMC & 41075.2 & 84 \\\hline
$n=10$ & ER & 16.8 & 130 \\
 & FMMC & $>10^6$ & 143 \\\hline
$n=20$ & ER & 19.20 & 315 \\
 & FMMC & $>10^6$ & 370 \\ \hline 
$n=25$ & ER & 20.00 & 403 \\
 & FMMC & $>10^6$ & 512\\
\hline
\end{tabular}
\caption{
FMMC vs ER on the small-world graph 
}
\label{tab:Comparison-SmallWorld}
\end{table}

\bc{In addition to FMMC, we consider Metropolis weights-based gossiping and fastest quantum gossiping (FQG) proposed by Jafarizadeh in \cite{jafarizadeh2020gossip}. In the Metropolis-based gossiping approach, each node $i$ wakes up \mg{with uniform probabilities} (i.e. $p_i^{m}=\frac{1}{n}$) and communicates with one of its neighbors $j \in \mathcal{N}_i \setminus \{i\}$ with probability $p_{j|i}^{m}= \frac{1}{\max\{d_i,d_j\}}$. The FQG, on the other hand, \mg{calculates} the wake-up ($p^f_i$) and conditional communication probabilities $(p^{f}_{j|i})$ at each agent $i\in \mathcal{N}$ by solving an SDP problem. This SDP is targeted to optimize the spectral gap of the expected iteration matrix. The method proposed in \cite{jafarizadeh2020gossip} for solving this SDP is a centralized algorithm; therefore, we made the comparisons among these methods in terms of the time required to compute the probability matrices $P^r$ and $P^{f}$ in a centralized manner. \mg{The entries of these matrices are computed according to $[P^f]_{ij}:=p^f_ip_{j|i}^f$ and $[P^r]_{ij}:=p^r_ip_{j|i}^r$. For the Metropolis weights, the probability matrix $P^m$ does not require any pre-computation time as it is only based on the degrees of the nodes and is assumed to be known. We also introduce the spectral gaps defined as} $\Delta_r=:1-\lambda_{n-1}(\overline{W}_{P^r})$, $\Delta_{f}:=1-\lambda_{n-1}(\overline{W}_{P^f})$, $\Delta_m:=1-\lambda_{n-1}(\overline{W}_{P^m})$ of the corresponding expected iteration matrices.}

\mg{In our next set of experiments, we consider barbell graphs and random graphs generated with the stochastic block model (SBM). The stochastic block model $\mbox{SBM}(n,k,p,q)$, also known as the planted partition model \cite{abbe2017community,lee2019review}, consist of $n$ nodes and $k$ clusters where each node in every cluster is connected to any other node in the same cluster with probability $p$, whereas the nodes that are not in the same cluster are connected with probability $q$.}

\begin{figure*}[ht]
\begin{footnotesize}
\centering
    \begin{subfigure}{0.3\textwidth}
    \centering
    \includegraphics[width=\textwidth]{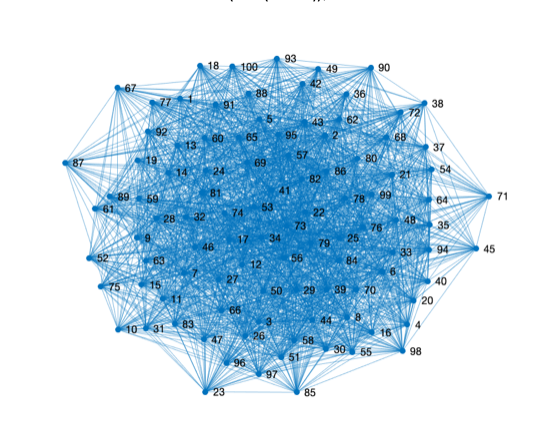}
    \caption{ {\small Small-world with $n$ nodes and $m$ links where $n=100$ and $\frac{2m}{(n^2-n)}\approx 0.02$.}}
    \label{fig: Smallworld}
    \end{subfigure}
    \begin{subfigure}{0.3\textwidth}
    \centering
    \includegraphics[width=\textwidth]{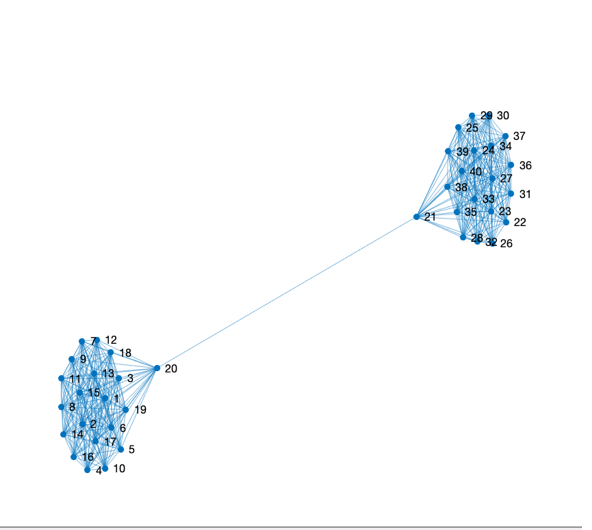}
    \caption{Barbell graph $K_{20}-K_{20}$.}
    \label{fig: barbell}
    \end{subfigure}
    \begin{subfigure}{0.3\textwidth}
    \centering
    \includegraphics[width=\textwidth]{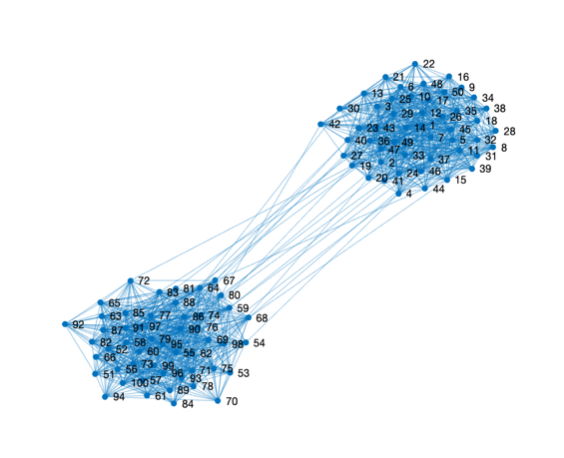}
    \caption{SBM(100,2,0.5,0.01)}
    \label{fig: SBM0_5_0_01}
    \end{subfigure}
    \begin{subfigure}{0.2\textwidth}
    \centering
    \includegraphics[width=\textwidth]{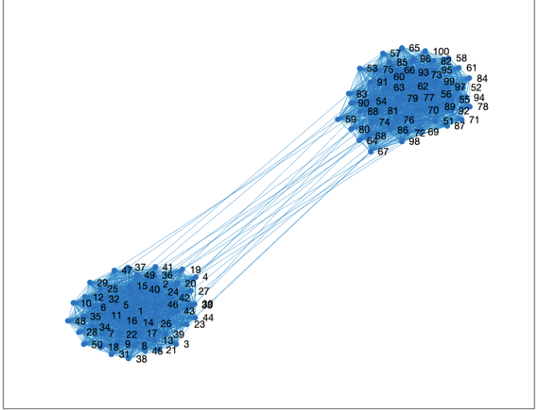}
    \caption{SBM(100,2,0.9,0.01)}
    \label{fig: SBM0_9_0_01}
    \end{subfigure}
    \begin{subfigure}{0.3\textwidth}
    \centering
    \includegraphics[width=\textwidth]{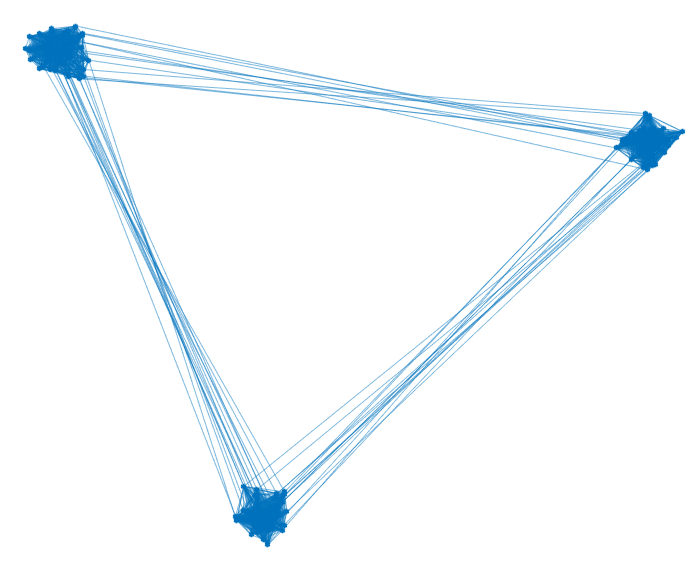}
    \caption{SBM(120,3,0.9,0.01)}
    \label{fig: SBM_3_0_01}
    \end{subfigure}
    \begin{subfigure}{0.3\textwidth}
    \centering
    \includegraphics[width=\textwidth]{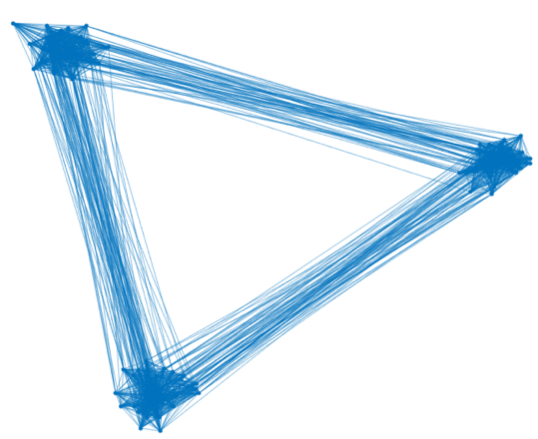}
    \caption{SBM(120,3,0.9,0.05)}
    \label{fig: SBM_3_0_05}
    \end{subfigure}
    \caption{\label{fig1} \bc{Different graph topologies used in numerical experiments.}}
    \end{footnotesize}
\end{figure*}

\bc{We summarize our results for barbell graph in Table \ref{Tab: SpecCompBarbell} and for $SBM(n,3,0.9,0.1)$ in Table \ref{Tab: SpecGapSBM3}. A gossiping algorithm will be faster if its spectral gap is larger. 
We observe that $\log(1/\Delta_{r})$ is \bc{smaller than $\log(1/\Delta_{m})$ and }larger than $1/\log(\Delta_{f})$ as $n$ is varied, therefore $\Delta_{r}$ is \bc{larger than $\Delta_{m}$} and is smaller than $\Delta_{f}$; so, we conclude that ER performs faster then Metropolis and slower than FQG. This is expected as effective-resistance weights are not \sa{optimized to increase the spectral gap} whereas FQG weights are targeted to optimize the spectral gap. However, when we look at the CPU time required to compute effective resistances and FQG weights, reported in the columns titled ``CPU Time ER" and ``CPU Time FQG", we observe that effective resistance weights can be computed faster \sa{as it only requires} a matrix inversion whereas FQG algorithm requires solving a semi-definite program (SDP).\footnote{We used the SDP solver \bc{SeDuMi} in the software \texttt{CVX} to compute FQG weights and matrix inversion function in Matlab to compute ER weights in a  centralized manner.}The advantage of \sa{using} the ER weights is that they are \sa{faster} to compute, and \sa{this effect becomes} more pronounced for larger graphs.
\sa{Moreover,} ER weights can be computed efficiently \sa{(with a {\it linear} convergence rate) and asynchronously} in the decentralized setting using the \sa{randomized} Kaczmarz algorithm. We note that solving SDPs in the decentralized setting is possible with subgradient methods as discussed in \cite{boyd2006randomized} but are typically much slower in the decentralized setting as they admit at most {\it sublinear} rates. \mg{This fact is also reflected in our results in Table \ref{Table: Barbell} and Table \ref{tab:Comparison-SmallWorld} for the FMMC method which used a subgradient method to compute the weights}. To summarize, we conclude that the Metropolis weights require no pre-computation but they are the slowest in terms of the spectral gap. ER is faster than Metropolis but slower than FQG weights; but the advantage is that computing ER weights require less CPU time. ER weights can also be computed efficiently in the decentralized setting with a linearly convergent algorithm.
}

\begin{table*}[ht]
\begin{footnotesize}
    \centering
    \begin{tabular}{|l|l|l|l|l|l|}
    \hline
         \textbf{Graph Topology} 
        & \textbf{\texttt{ER-Kac}} 
        & \textbf{\texttt{ER-Ex}}
        & \textbf{\texttt{Uniform}}
        & \textbf{\texttt{Metropolis}}
       \\
        \hline
       \textbf{Small-world, $n=100$}
        & \textbf{530 $(\pm 28.92)$}
        & \textbf{527 $(\pm 29.29)$}
        & 801 $(\pm 219.22)$
        & 623 $(\pm 43.10)$ 
        \\ 
        \hline 
        \textbf{$\mathbf{K_{20}-K_{20}}$}
        & \textbf{1577} $(\pm 563.62)$
        & \textbf{1525} $(\pm 505.93)$ 
        & 8661 $(\pm 3460.71)$ 
        & 8341 $(\pm 3328.85)$
        \\  
        \hline
        \textbf{SBM(100,2,0.5,0.01)}
        & \textbf{749 $(\pm 306.76)$} 
        & \textbf{721 $(\pm 268.02)$}
        & 1213 $(\pm 363.60)$
        & 914 $(\pm 422.66)$ 
        \\ 
        \hline
        \textbf{SBM(100,2,0.9,0.01)} 
        & \textbf{835 $(\pm 436.50)$}
        & \textbf{814 $(\pm 439.41)$} 
        & 1459 $(\pm 643.22)$
        & 1248 $(\pm 960.33)$ 
        \\
        \hline
        \textbf{SBM(120,3,0.9,0.01)}
        & \textbf{1022 $(\pm 383.75)$}
        & \textbf{1003 $(\pm 349.30)$}
        & 2079 $(\pm 604.95)$
        & 1541 $(\pm 762.64)$
        \\ 
        \hline
        \textbf{SBM(120,3,0.9,0.05)}
        & \textbf{646} $(\pm 59.98)$
        & \textbf{639} $(\pm 65.43)$
        & 1185 $(\pm 158.95)$
        & 674 $(\pm 69.33)$
        \\
        \hline
    \end{tabular}
    \caption{ \response{The comparison of mean and the standard deviation of the wake-ups required for ER, Metropolis, and uniform asynchronous gossiping algorithms based on 250 runs.}}
    \label{Table: WakeUp Comparison}
    \end{footnotesize}
    
\end{table*}

\begin{table}[ht!]
\centering
\begin{small}
\resizebox{0.8\columnwidth}{!}{
\begin{tabular}{|c|c|c|c|c|c|}
\hline
\textbf{$\tilde{n}$} & $\log(1/\Delta_{r})$ & $\log(1/\Delta_{m})$ & $\log(1/\Delta_{f})$ & \textbf{\begin{tabular}[c]{@{}l@{}l@{}l@{}} CPU Time \\ ER (in secs)\end{tabular}}
& \textbf{\begin{tabular}[c]{@{}l@{}l@{}l@{}} CPU Time \\ FQG (in secs)\end{tabular}} \\ \hline
20 & 5.914 & 7.076  & 4.407 & $\leq 0.01$ & 2.50   \\ \hline
24 & 6.292 & 7.599  & 4.605 & $\leq 0.01$ & 3.10   \\ \hline
28 & 6.610 & 8.043  & 4.771 & $\leq 0.01$ & 3.48   \\ \hline
32 & 6.884 & 8.429  & 4.913 & $\leq 0.01$ & 6.47   \\ \hline
36 & 7.125 & 8.771  & 5.037 & $\leq 0.01$ & 8.63   \\ \hline
40 & 7.340 & 9.078  & 5.147 & $\leq 0.01$ & 12.33  \\ \hline
44 & 7.534 & 9.357  & 5.247 & $\leq 0.01$ & 13.45  \\ \hline
52 & 7.873 & 9.846  & 5.421 & $\leq 0.01$ & 21.79  \\ \hline
72 & 8.532 & 10.803 & 5.756 & $\leq 0.01$ & 106.22 \\ \hline
\end{tabular}
}
\caption{\response{Comparison of spectral gaps  $\Delta_r, \Delta_m$, and $\Delta_f$ and CPU times (in seconds) required to compute the communication $(p_{i|j})$ and wake-up probabilities $(p_i)$ on barbell graphs $K_{\tilde{n}}-K_{\tilde{n}}$.}}
\label{Tab: SpecCompBarbell}
\end{small}
\end{table}

\begin{table}[ht!]
\centering
\begin{small}
\resizebox{0.80\columnwidth}{!}{
\begin{tabular}{|c|c|c|c|c|c|}
\hline
\textbf{n} & $\log(1/\Delta_{r})$ & $\log(1/\Delta_{m})$ & $\log(1/\Delta_{f})$ & 
\textbf{\begin{tabular}[c]{@{}l@{}l@{}l@{}} CPU Time \\ ER (in secs)\end{tabular}}
& \textbf{\begin{tabular}[c]{@{}l@{}l@{}l@{}} CPU Time \\ FQG (in secs)\end{tabular}} \\ \hline
22 & 6.490 & 7.150 & 5.457 & $\leq 0.01$ & 3.02  \\ \hline
24 & 6.276 & 7.356 & 5.072 & $\leq 0.01$ & 4.09  \\ \hline
30 & 6.743 & 8.049 & 5.236 & $\leq 0.01$ & 4.02  \\ \hline
32 & 7.104 & 8.328 & 5.573 & $\leq 0.01$ & 6.98  \\ \hline
36 & 6.498 & 7.524 & 5.302 & $\leq 0.01$ & 6.74  \\ \hline
38 & 6.268 & 6.959 & 4.974 & $\leq 0.01$ & 8.58  \\ \hline
40 & 6.613 & 7.475 & 5.264 & $\leq 0.01$ & 8.97  \\ \hline
44 & 6.621 & 7.420 & 5.265 & $\leq 0.01$ & 11.04 \\ \hline
48 & 6.904 & 7.814 & 5.361 & $\leq 0.01$ & 13.58 \\ \hline
50 & 6.969 & 7.947 & 5.505 & $\leq 0.01$ & 15.05 \\ \hline
52 & 7.166 & 8.249 & 5.566 & $\leq 0.01$ & 16.93 \\ \hline
54 & 6.934 & 7.840 & 5.474 & $\leq 0.01$ & 21.39 \\ \hline
\end{tabular}
}
\caption{\response{The comparison of spectral gaps $\Delta_r, \Delta_m$, and $\Delta_f$ of the iteration matrices and CPU times (in secs) required to compute $P^r$ and $P^f$ on $SBM(n,3,0.9,0.01)$.}}
\label{Tab: SpecGapSBM3}
\end{small}
\end{table}

\bc{Lastly, we compare the performance of ER-based asynchronous gossiping with Metropolis weights-based asynchronous gossiping (\texttt{Metropolis}) and classical asynchronous gossiping (\texttt{Uniform}) on small-world, barbell, and random graphs generated with the stochastic block model, $\mbox{SBM}(n,k,p,q)$. We considered two types of ER-based gossiping algorithms: (i) The first algorithm 
\sa{\texttt{ER-Ex}} uses the exact effective resistance probabilities that are computed based on calculating the pseudo-inverse of the Laplacian matrix (with a centralized approach based on standard matrix inversion techniques), (ii) The second algorithm \sa{\texttt{ER-Kac}} is based on the effective-resistance weights approximated by the decentralized Kacmarz method.
}
\bc{
In the experiments, each node $i$ possesses an initial vector $y_i^{(0)} \in \mathbb{R}^5$ and the goal is to approximate the node averages $\bar{y} = \frac{1}{n} \sum_{i=1}^n y_i^{(0)}$. We draw the data $y_{i}^{(0)}$ randomly according to the standard multi-variate normal distribution admitting a zero mean and a unit covariance matrix. In each trial, we record the number of wake-ups required to obtain the relative accuracy $\frac{\sum_{i=1}^{n}\Vert y_i^{(k)}-\bar{y}\Vert^2}{n\Vert \bar{y}\Vert^2}\leq \varepsilon$. We set $\varepsilon = 5\times 10^{-6}$ and generated 250 independent runs. We calculated the average and the standard deviation of the wake-ups among these 250 runs on SBM, barbell graph, and small-world graphs. We presented our results in the Table \ref{Table: WakeUp Comparison} (see Figure \ref{fig1} for the details of these graphs).
}

\bc{We observe that effective-resistance based algorithms (\texttt{ER-Kac} and \texttt{ER-Ex}) improve clearly upon the uniform (uniform weights-based) gossiping on all of the graph types in terms of both average wake-ups required and the standard deviation of the wake-ups. \sa{When we compare \texttt{ER-Kac} and \texttt{ER-Ex} with Metropolis weights, we observe that \texttt{ER-Kac} and \texttt{ER-Ex} are more efficient compared to Metropolis weights in the sense that they require a smaller number of wake-ups on average with a smaller standard deviation.} The improvement is more pronounced for the graphs in Fig. \ref{fig: barbell}--\ref{fig: SBM_3_0_01}. 
These experimental results illustrate the effectiveness of our approach compared to existing approaches on a number of \sa{random} graph topologies that can arise in practice.}

\subsection{Effective resistance-based DPGA-W and EXTRA}
\begin{figure}[h!]
\begin{center}
	\includegraphics[width=0.9\linewidth]{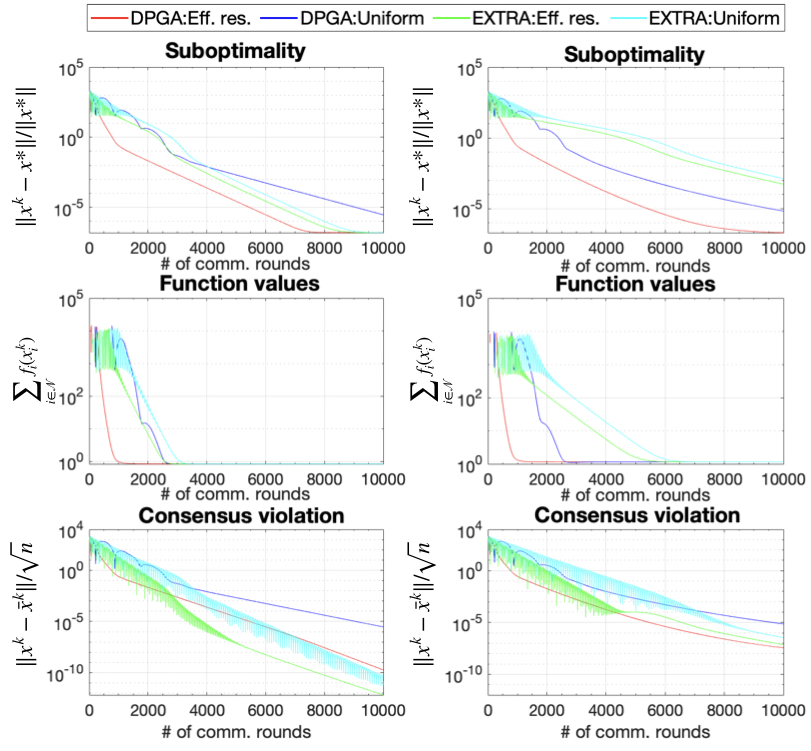}
\caption{\small The suboptimality, function value, and the difference from average versus \bc{the number of communication rounds}, \mg{based on} logistic regression using DPGA-W and EXTRA algorithms with resistance weights and uniform probability weights on barbell graph $K_{10}-K_{10}$. \textbf{Left:} Data of the logistic regression model is sampled using $\sigma=1$, \textbf{Right:}  Data is sampled using $\sigma=2$.}
\label{fig:dpga_n20}
\end{center}
\end{figure}

\begin{figure}[h!]
\begin{center}
	\includegraphics[width=0.9\linewidth]{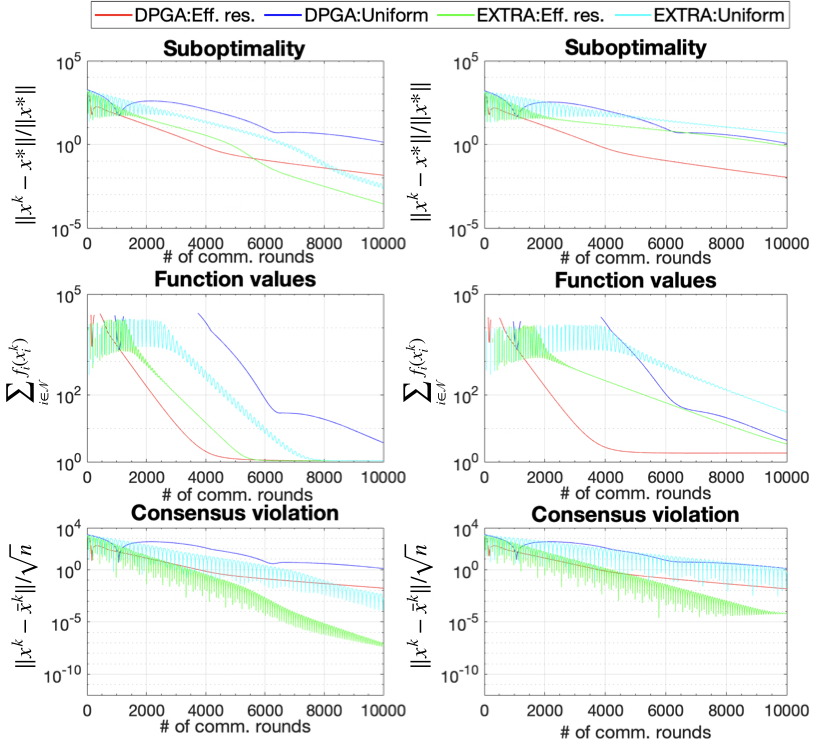}
\caption{\small The suboptimality, function value, and difference from average \bc{versus average number of communication rounds}, \bc{based on} logistic regression using DPGA-W and EXTRA algorithms with resistance weights and uniform probability weights on the barbell graph $K_{20}-K_{20}$. \textbf{Left:} Data of the logistic regression model is sampled using $\sigma=1$, \textbf{Right:}  Data is sampled using $\sigma=2$.}
\label{fig:dpga_n40}
\end{center}
\vspace*{-8mm}
\end{figure}

We implemented our 
ER-based communication framework into the state of the art distributed algorithms: DPGA-W~\cite{aybat2018distributed} and EXTRA~\cite{shi2015extra} 
to solve regularized logistic regression problems over a barbell graph $K_{\tilde{n}}-K_{\tilde{n}}$ with $n=2\tilde{n}$ nodes: We minimize $\min_{x\in\reals^p}\sum_{i=1}^n f_i(x)$ with
{\small
\begin{align}
    \label{eq:log_reg}
    f_i(x)\triangleq\tfrac{1}{2n}\norm{x}^2+\tfrac{1}{N_s}\sum_{\ell=1}^{N_s}\log(1+\exp^{-b_{i\ell} a_{i\ell}^\top x}),
\end{align}}%
where $N_s$ is the number of samples at each node, $\{(a_{i\ell},b_{i\ell})\}_{\ell=1}^{N_s}\subset\reals^p\times\{-1,1\}$ for $i\in\cN$ denote the set of feature vectors and corresponding labels. {We let $p=20$ and $N_s=5$. For each $n\in\{20,40\}$ and $\sigma\in\{1,2\}$, we randomly generated $20$ i.i.d. instances of the problem in~\eqref{eq:log_reg} by sampling
$a_{i\ell}\sim \textbf{N}(\one,\sigma^2\bI)$ independently from the normal distribution and setting $b_{i\ell}=-1$ if $1/(1+e^{-a_{i\ell}^\top\one})\leq 0.55$ and to $+1$ otherwise}. 
Both algorithms are terminated after $10^4$ iterations. For benchmark, we also solved each instance of \eqref{eq:log_reg} using MOSEK~\cite{mosek} within CVX~\cite{cvx}. We initialized the iterates uniformly sampling each $p$ components from the $[500, 510]$ interval for nodes in one $K_\tn$, and from $[-500, -490]$ for nodes in the the other $K_\tn$. The results for $n=20$ and $n=40$ are displayed in Fig.~\ref{fig:dpga_n20} and Fig.~\ref{fig:dpga_n40}, respectively. We plotted relative suboptimality $\norm{\bx^k-\bx^*}/\norm{\bx^*}$, function value sequence $\sum_{i\in\cN}f_i(x_i^k)$ for the range $[0,~10^5]$, and consensus violation $\norm{\bx^k-\bar{\bx}^k}/\sqrt{n}$, where $k$ denotes the (synchronous) communication round counter -- in each communication round neighboring nodes communicate among each other synchronously once -- and $\bx^k=[x_i^k]_{i\in\cN}$ denotes the $k$th iterate; 
moreover, $\bar{\bx}^k=\one\otimes\bar{x}^k$, $\bar{x}^k=\sum_{i\in\cN}x_i^k/n$, $\bx^*\triangleq\one\otimes x^*$ and $x^*$ is the minimizer to \eqref{eq:log_reg}. 

Both DPGA-W\footnote{In DPGA-W stepsize parameter $\gamma_i$ is set to $1/\norm{\omega_i}$ for $i\in\cN$ -- see~\cite{aybat2018distributed}.} and EXTRA uses a communication matrix $W$ that encodes the network topology. DPGA-W uses node-specific step-sizes initialized at $\approx 1/L_i$ for $i\in\cN$, where $L_i$ denotes the Lipschitz constant of $\grad f_i$, we adopted the adaptive step-size strategy described in \cite[Sec. III.D]{aybat2018distributed}; and for EXTRA, we choose the constant step-size, common for all nodes, as suggested in \cite{shi2015extra}, i.e., we choose the step size as  $2\lambda_{\min}(\tilde{W})/\max_{i\in\cN}{L_i}$, where $\tilde{W}=(\bI+W)/2$. 

{For both algorithms, we compared two choices of $W$: $W^u$ based on uniform edge weights, and $W^r$ based on effective resistances. In DPGA-W, the graph Laplacian is adopted for uniform weights, i.e., $W^u = W^{u,\text{\tiny DPGA-W}}\triangleq\cL$, while for the 
ER-based weights, we set $W^r = W^{r,\text{\tiny DPGA-W}}$ where $W^{r,\text{\tiny DPGA-W}}_{ii}\triangleq\sum_{j\in\cN_i}R_{ij}$ for $i\in\cN$ and $W^{r,\text{\tiny DPGA-W}}_{ij}=-R_{ij}$ for $(i,j)\in\cE$ and $0$ otherwise. For EXTRA, $W^{u,\text{\tiny EXTRA}}=\bI-\cL/\tau$ where $\tau=\lambda_{\max}(\cL)/2+\varepsilon$ where $\lambda_{\max}$ denotes the largest eigenvalue}; on the other hand, $W^{r,\text{\tiny EXTRA}}=\bI-W^{r,\text{\tiny DPGA-W}}/\tau$ where $\tau=\lambda_{\max}(W^{r,\text{\tiny DPGA-W}})/2+\varepsilon$ for $\varepsilon=0.01$.

{Figures \ref{fig:dpga_n20} and \ref{fig:dpga_n40} illustrate the performance comparison of both DPGA-W and EXTRA algorithms with effective resistance and uniform weights in terms of suboptimality, convergence in function values and consensus violation for the barbell graph $K_{10}-K_{10}$ and $K_{20}-K_{20}$ respectively -- the reported results are averages over the 20 problem instances. The 
subfigures on the left of Figures \ref{fig:dpga_n20} and \ref{fig:dpga_n40} are for noise level $\sigma=1$ whereas 
those on the right are for 
$\sigma=2$. In Figures \ref{fig:dpga_n20} and \ref{fig:dpga_n40}, we observe that 
using ER weights improves upon the uniform weights for both EXTRA and DPGA-W methods consistently to solve the logistic regression problem in terms of suboptimality, function values and consensus violation significantly. We also observe that with noisier data, DPGA-W works typically faster than EXTRA in terms of function values and suboptimality. This is because when noise level $\sigma$ gets larger, the local Lipschitz constant $L_i$ of the nodes demonstrate higher variability, and DPGA-W adapts to this variability as it uses a step size that is different at each node in a way to adapt to $L_i$, whereas EXTRA uses a constant step size that is the same for 
all nodes. On the other hand, in terms of consensus violation, we see that EXTRA with 
ER weights typically outperforms DPGA-W with 
ER weights. 
}
\section{Conclusions}\label{sec-future}
We obtained a number of theoretical guarantees for ER gossiping algorithms for the consensus problem for $c$-barbell graphs \bc{and} barbell graphs, and for arbitrary graphs depending on their \response{diameter}. The results fill a gap between the theory and practice of these methods. We also showed that these methods are effective for solving the consensus problem in practice over barbell graphs and small-world graphs.
We provided numerical experiments 
demonstrating that
using ER gossiping within EXTRA and DPGA-W methods improves their practical performance in terms of communication efficiency.


\section*{Acknowledgments}
Bugra Can and Mert G\"urb\"uzbalaban's research were supported by the the Office of Naval Research Award Number N00014-21-1-2244, and the grants National Science Foundation (NSF) CCF-1814888 and NSF DMS-2053485. N. Serhat Aybat's research is supported by the grants NSF CMMI-1635106 and ARO W911NF-17-1-0298.

\appendix

\vspace{-1mm}
\section{Proof of Propositions \ref{prop: c-barbell graph} and \ref{prop: eig-order-imprv}}\label{sec: Prop 6and9}
\begin{proof}[Proof of Proposition \ref{prop: c-barbell graph}] 
The proof is 
\bc{based on finding the subset $S$ of the vertex set of $c$-barbell graph \mg{that determines the conductance}, i.e. that solves the minimization problem \eqref{eq: conductance}. First, for any given $\cG=(\cN,\cE)$}, the conductance of a subset $S\subset \cN$ with respect to the probability transition matrix $W$ is defined as
\begin{equation}\label{cond: set_conduc}
\Phi_S (W) \triangleq \frac{1}{\pi(S)}\sum_{i\in S, j\in S^C}\pi_{i}W_{ij}. 
\end{equation} 
\bc{Notice that the definition \eqref{eq: conductance} implies that we have $\Phi(W)=\min_{\sa{S\subset\cN:}\pi(S)\in (0,1/2]}\Phi_S(W)$.}\footnote{\response{This follows after straightforward computations based on the the fact that the Markov chain with transition matrix $W$ and stationary distribution $\pi$ is a reversible Markov chain, i.e. $\pi(S)\Phi_S(W) = \pi(S^c)\Phi_{S^c}(W)$ for any $S$ with $\pi(S)\in (0,1)$.}}
With slight abuse of notation, for  a subgraph $\mathcal{H}_0$ with a vertex set $S_0$, we define $\Phi_{\mathcal{H}_0} (W) \triangleq \Phi_{S_0} (W)$. 
We say that a vertex set $S\subset \mathcal{N}$ on graph  $\mathcal{G}=(\mathcal{N},\mathcal{E},w)$ is a \textit{one-cut set} if its complement $\mathcal{N}\setminus S$ is a connected subgraph of $\mathcal{G}$. Similarly, we define \textit{two-cut set} $S_2\subset  \mathcal{N}$ to be a set whose complement $\mathcal{N}\setminus S_2$ consists of two disjoint  non-empty connected subgraphs $\mathcal{H}_1$ and $\mathcal{H}_2$ of $\mathcal{G}$. We define
\begin{small}
\begin{equation}\label{def-G1}\mathcal{G}_{1} \triangleq \text{the left-most clique of the }c\text{-barbell graph}.
\end{equation}
\end{small}
For $c_0\in [2,c]$, we also define 
\begin{small}
\begin{align}
\mathcal{G}_{c_0} \triangleq & ~\mbox{$c_0$-barbell subgraph that includes the left-most}\nonumber\\
&~\mbox{$c_0$ cliques of the $c$-barbell graph.
}\label{def: subgraph_G_c}
\end{align} 
\end{small}
\bc{Note that} matrices $W_{P^u}$ and $W_{P^r}$ are symmetric and Markov chains with these transition matrices have the uniform distribution as a stationary distribution.
\bc{Therefore, Lemmas \ref{Prop: one cut} and \ref{Prop: one cut at edge} provided in Appendix~\ref{sec: Supporting Results} imply that a set $S$ with minimal conductance should be a one-cut set and has to be  given by the vertices of a subgraph $\cG_{c_0}$ for some $c_0 \in [1,c]$ for both $W_{P^u}$ and $W_{P^r}$.}
\bc{The conductance of one-cut subgraphs with respect to these transition matrices can be computed explicitly (see Proof of Lemma \ref{Prop: one cut at edge} for details):}
\begin{align} \label{cond: C-Kn uni and res}
\Phi_{\cG_{c_0}}(\overline{W}_{P^u})=\frac{1}{c_0}\frac{1}{c\tn^3},\;\;
\Phi_{\cG_{c_0}}(\overline{W}_{P^r})=\frac{1}{2c_0\tn(c\tn-1)}.
\end{align} 
Both of the expressions at \eqref{cond: C-Kn uni and res} are minimized for the choice of $c_0=\lfloor \frac{c}{2} \rfloor $. Therefore, the minimal conductance is attained for the subgraph
$\cG_{\lfloor \frac{c}{2}\rfloor}$. 
Plugging $c_0 =\lfloor \frac{c}{2}\rfloor$ into the expressions above yields the graph conductance values at \eqref{Cond: c-barbell}. The bounds \eqref{bound: c_barbell_uni} and \eqref{bound:c_barbell_res} follow from Theorem \ref{thm-gossip-compare} and inequalities \eqref{eq:lambda_bound}.
\end{proof}
\begin{proof}[Proof of Proposition \ref{prop: eig-order-imprv}]
It follows from Corollary \ref{cor: stochastic e-values} and Lemma \ref{lemma: Values} in Appendix~\ref{sec: Supporting Results} that the second largest eigenvalues of 
$\bar{W}_{P^u}$ and $\bar{W}_{P^r}$ are given by:
$\lambda_{n-1}(\overline{W}_{P^{u}})= 1-\frac{8}{n^2(n-2)}+\Theta(\frac{1}{n^4})$ and 
	$\lambda_{n-1}(\overline{W}_{P^{r}})= 1- \frac{1}{n(n-1)}-\Theta (\frac{1}{n^3})$.
This implies directly 
   $\lambda_{n-1}(\overline{W}_{P^r}) = 1 - \Theta(\frac{1}{n^2})$ and 
   $\lambda_{n-1}(\overline{W}_{P^u}) = 1 -\Theta(\frac{1}{n^3})$,
which completes the proof.
\end{proof}
\vspace{-3mm}
\section{Supporting Results}\label{sec: Supporting Results}
\begin{lemma}\label{Prop: one cut} Consider a reversible Markov chain on a $c$-barbell graph with a uniform stationary distribution. Let
$\mathcal{H}_0$ be a subgraph of $\mathcal{G}$ whose vertex set is a non-empty two-cut set
$\mathcal{S}_0$ satisfying $|\mathcal{S}_0|\leq \frac{|\mathcal{N}|}{2}$. 
Then, there exists another subgraph $\tilde{\mathcal{H}}_0$ of $\cG$ such that $\Phi_{\tilde{\mathcal{H}}_0}(W)<\Phi_{\mathcal{H}_0}(W)$. 
\end{lemma} 
\begin{proof} 
Let $C_1$ and $C_2$ be the vertex sets of two disjoint non-empty connected subgraphs within $\cN\setminus \cS_0$ satisfying $\mathcal{N}=C_1 \cup \mathcal{S}_0 \cup C_2$. Note that $C_1 \cap C_2 = \emptyset$ implies either $|C_1 \cup \mathcal{S}_0|\leq \frac{|\mathcal{N}|}{2}$ or $\vert C_2 \vert \leq \frac{\vert \mathcal{N}\vert}{2}$. Using the fact that the transition matrix $W$ of a reversible Markov chain with a uniform stationary distribution is symmetric, the definition \eqref{cond: set_conduc} implies $\Phi_{C_1 \cup \mathcal{S}_0}(W)=\Phi_{C_2 }(W)$.
Without loss of generality, choose $\tilde{\mathcal{H}}_0$ to be the subgraph with vertices $\tilde{\mathcal{S}_0}=C_1\cup \mathcal{S}_0$ with $\vert C_1 \cup  \mathcal{S}_0\vert \leq \frac{ \vert \mathcal{N} \vert}{2}$ (otherwise, pick the subgraph with vertex set $C_2$ instead), then
\begin{small}
\begin{align*}
    \Phi_{\mathcal{H}_0}(W)&=\frac{1}{|\mathcal{S}_0|}\Big(\sum_{\substack{i\in \mathcal{S}_0\\ j\in C_1}} W_{ij}+\sum_{\substack{i\in \mathcal{S}_0\\ j\in C_2}}W_{ij} \Big)\\
    &>  \frac{1}{|\mathcal{S}_0|}\sum_{\substack{i\in \mathcal{S}_0\\ j\in C_2}}W_{ij}
    >  \frac{1}{|\tilde{\mathcal{S}}_0|}\sum_{\substack{i\in\tilde{\mathcal{S}}_0\\j\in C_2}} W_{ij}
    =\Phi_{\tilde{\mathcal{H}}_0}(W),
\end{align*}
\end{small}
which proves Lemma \ref{Prop: one cut}.
\end{proof} 
\begin{lemma} \label{Prop: one cut at edge} Consider a Markov chain on a $c$-barbell graph with a probability transition matrix $W$. If $W=\overline{W}_{P^u}$ or $W=\overline{W}_{P^r}$, then for any subgraph $\mathcal{H}_0$ having a one-cut vertex set $\mathcal{S}_0$, there exists a subgraph $\mathcal{G}_{c_0}$ for some $c_0 \in [1,c]$ such that $\Phi_{\mathcal{G}_{c_0}}(W) \leq \Phi_{\mathcal{H}_0}(W)$ where $\mathcal{G}_{c_0}$ is defined by \eqref{def-G1} and \eqref{def: subgraph_G_c}. 
\end{lemma} 
\begin{proof} For any subgraph $\mathcal{H}_0$ having a one-cut vertex set $\mathcal{S}_0$, we can always a find a subgraph $\mathcal{G}_{c_0}$ with vertex set $\mathcal{V}_{c_0}$ for some $c_0\in [1,c]$ such that either $\mathcal{V}_{c_0-1} \subset \mathcal{S}_0 \subset \mathcal{V}_{c_0}$ or $\mathcal{V}_{c_0-1}\subset S^c_0\subset \mathcal{V}_{c_0}$ (with the convention that $\mathcal{G}_{c_0}$ is a singleton graph with a vertex set $\mathcal{V}_{0}$ consisting of a single node). 
Let $\mathcal{H}_{0}^c$ be the subgraph with vertex set $S^c_0$. Since $\Phi_{\mathcal{H}_0}(W)=\Phi_{\mathcal{H}^c_0}(W)$ for both $W = \overline{W}_{P^r}$ and $W=\overline{W}_{P^u}$, without loss of generality we can assume that $\mathcal{H}_{0}$ satisfies the property $\mathcal{V}_{c_0-1} \subset \mathcal{S}_0 \subset \mathcal{V}_{c_0}$ (otherwise, we can replace $\mathcal{H}_0$ with $\mathcal{H}_0^c$ in the proof below).
It follows after a straightforward computation (similar to the proof technique of Lemma \ref{lemma: Values}) that transition probability matrices $\overline{W}_{P^u}$ and $\overline{W}_{P^r}$ on $c-K_\tn$ admit the explicit formula
    $[\overline{W}_{P^u}]_{i^*j^*}= \frac{1}{c\tn^2}$, $[\overline{W}_{P^u}]_{i^*j}=\frac{1}{2c\tn^2}\Big(\frac{2\tn-1}{\tn-1} \Big)$, $[\overline{W}_{P^u}]_{ij}=\frac{1}{c\tn(\tn-1)}$,
whereas 
    $[\overline{W}_{P^r}]_{i^*j^*}=\frac{1}{2(c\tn-1)}$, $[\overline{W}_{P^r}]_{i^*j}=\frac{1}{\tn(c\tn-1)}$, $[\overline{W}_{P^r}]_{ij}=\frac{1}{\tn(c\tn-1)}$,
where $i^*$ and $j^*$ denote two adjacent nodes belonging to different complete subgraphs of $c-K_{\tn}$, i.e., those with degree $\tn$, and $(i,j)\in\cE$ or $(i^*,j)\in\cE$ such that $i$ and $j$ denote nodes in $c-K_{\tn}$ with degree $\tn-1$.
Note $[W_{P^r}]_{i^*j^*}$ is greater than $[W_{P^u}]_{i^*j^*}$ as in the $K_\tn-K_\tn$ case. Hence, for $W=\overline{W}_{P^u}$,\vspace*{-2mm} 
\begin{small}
\begin{eqnarray*} 
\Phi_{\mathcal{H}_0}(\overline{W}_{P^u})= \frac{1}{|\mathcal{S}_0|}\sum_{\substack{i\in \mathcal{S}_0\\j\in \mathcal{S}_0^c }}[\overline{W}_{P^u}]_{ij} > \frac{1}{c_0\tn} \frac{1}{c\tn^2}
= \Phi_{\mathcal{G}_{c_0}}(\overline{W}_{P^u}).
\end{eqnarray*}
\end{small}
\vspace*{-1mm}
In the case of  $W=\overline{W}_{P^r}$, let $\cP_0\subset\cS_0$ be the subset of nodes in the subgraph $K_{\tn}$ that contains nodes from both $\cS_0$ and $\cS_0^C$ -- if no such $K_{\tn}$ exists, then $\cS_0$ corresponds to a subgraph $\mathcal{G}_{c_0}$ for some $c_0 \in [1,c]$. Now consider the former case, let us denote $m_0 \triangleq |\mathcal{P}_0|<\tn$. The number of edges between $\mathcal{P}_0$ and $\mathcal{S}_0^C$ is given by $m_0(\tn-m_0)$. This is due to the fact that each 
node in $\cP_0$ has exactly $(\tn-m_0)$ many edges that connect $\mathcal{S}_0$ to its complement. We have also $m_0(\tn-m_0)\geq \frac{\tn}{2}$ for $\tn\geq 2$. This yields 
$\Phi_{\mathcal{H}_0}(\overline{W}_{P^r})=\frac{1}{|\mathcal{S}_0|}\sum_{\substack{i\in \mathcal{S}_0\\ j\in \mathcal{S}_0^c}} [\overline{W}
_{P^r}]_{ij}
            \geq \frac{1}{|\mathcal{S}_0|} \frac{m_0(\tn-m_0)}{\tn(c\tn-1)}\\
            \geq \frac{1}{|c_0\tn|} \frac{1}{2(c\tn-1)}
         =\Phi_{\mathcal{G}_{c_0}}(W_{P^r})$.
\end{proof}
\begin{corollary}\label{cor: stochastic e-values} Under the setting of Proposition \ref{prop: gen. eigenvalue}, assume that the weight matrix $w$ is normalized, i.e., $\sum_{j=1}^{n}w_{ij}=1$ for all $i\in \mathcal{N}$. Then $W=w$ is a doubly stochastic matrix and the eigenvalues of $W$ become
\begin{itemize}
    \item $\lambda_{a}=1$ with multiplicity one,
    \item $\lambda_{b}=-1+(A+G)+F$ with multiplicity one,
    \item $\lambda_{c}=D-C$ with multiplicity $2\tn-4$,
    \item $\lambda_{\pm}=\frac{1}{2}\Big(F+G-A\,\pm\,\sqrt{S} \Big)$,
\end{itemize}
where $A,B,C,D,E,F,G$ and $S$ are as in Proposition \ref{prop: gen. eigenvalue}. Moreover, $\lambda_+$ satisfies\vspace*{-1mm}
{\small
\begin{equation}\label{eq-lambda-star}
\lambda_{+}=\frac{1}{2}\Big(F+G-A\,+\,\sqrt{(F-G+A)^2+4BE} \Big),
\end{equation}}%
\vspace*{-1mm}
\noindent and is the second largest eigenvalue, i.e. 
$\lambda_{n-1}(W) = \lambda_+$.
\end{corollary}
\begin{proof} 
Since $w$ is normalized, Proposition \ref{prop: gen. eigenvalue} applies with  $A+G+E=1$ and $B+F=1$. Thus eigenvalues simplify to the forms given in the statement. Note that
    $\sqrt{S}=\sqrt{(F+G-A)^2-4(FG-BE-AF)} 
            =\sqrt{(F-G+A)^2+4BE}\geq 0$.
Therefore, $\lambda_+$ satisfies 
\eqref{eq-lambda-star}. 
$\lambda_a = 1$ is the unique largest eigenvalue since $W$ is stochastic. It remains to show that $\lambda_+$ is the second largest eigenvalue. Using \eqref{eq-lambda-star}, we can write
$
\lambda_{+}\geq \frac{1}{2}\big(F+G-A+ \vert F-G+A\vert \big).
$
There are two cases: $F \geq (G-A)$ or $F < (G-A)$. 
In both cases, we observe 
$\lambda_{+}\geq F\geq 0$. Since $A+G+E=1$, we also have $A+G-1 = -E \leq 0$. Therefore $\lambda_b = F-E \leq F \leq  \lambda_+.$
Furthermore, $\lambda_c = D-C\leq F=D + (\tn-2)C$ since $C\geq 0$; therefore $\lambda_c \leq F \leq \lambda_+$. Finally, $\lambda_+ \geq 0$ since $S\geq 0$. 
Thus, $\lambda_{+}$ is non-negative and is the second largest eigenvalue.
\end{proof}
\begin{lemma}\label{lemma: Values} Consider the setting of Proposition \ref{prop: gen. eigenvalue}:
\begin{itemize}
\item [$(i)$] If $W = \overline{W}_{P^u}$, then Proposition \ref{prop: gen. eigenvalue} applies with $A = A^u$, $B = B^u$, $C = C^u$, $D = D^u$ and $G=G^u$ where 
\begin{align*}
A^{u}&= \frac{2}{n^2},~ 
B^{u}=\frac{n-1}{n^2(0.5n-1)},~ 
C^{u}= \frac{2}{n(n-2)},\\
D^{u}&= \frac{n^3-3n^2+2n+2}{n^2(n-2)},~
G^u =1-\frac{n+1}{n^2}. 
\end{align*}
The second largest eigenvalue of $\overline{W}_{P^u}$ is given by
	$\lambda_{n-1}(\bar{W}_{P^{u}})=1-\frac{n^2+n-8}{2n^2(n-2)}+\frac{1}{8}\sqrt{S_n^u}
	=  1- \frac{8}{n^2(n-2)}+\Theta(\frac{1}{n^4})$,
for $S_n^u=\frac{4n^3+24n^2-156n+192}{(0.5n-1)^2n^3}$. 
\item [$(ii)$] If $W = \overline{W}_{P^r}$, then Proposition \ref{prop: gen. eigenvalue} applies with $A = A^r$, $B = B^r$, $C = C^r$, $D = D^r$ and $G=G^r$ where
\begin{align*}
A^{r}&= \frac{1}{2(n-1)},~
 B^{r}=\frac{2}{n(n-1)},~ C^{r}= \frac{2}{n(n-1)},\\
D^{r}&= \frac{n^2-2n+2}{n(n-1)},~ G^r=1-\frac{1.5n-2}{n(n-1)}.
\end{align*}
Moreover, the second largest eigenvalue of $\overline{W}_{P^u}$ is given by
		$\lambda_{n-1}(\bar{W}_{P^{r}})=1-\frac{1}{(n-1)}+\frac{1}{2}\sqrt{S_n^r}
		= 1- \frac{1}{n(n-1)}-\Theta (\frac{1}{n^3})$,
for $S_n^r=\frac{4n-8}{n(n-1)^2}$.
\end{itemize}
\end{lemma}
\vspace{-2mm}
\begin{proof}[Proof of Lemma \ref{lemma: Values}]  We first compute the entries of both $P^u$ and $P^r$ matrices explicitly for the barbell graph (i.e. $K_\tn-K_\tn$). Former one can be found directly from degrees of the nodes: $P^{u}_{ij}=\frac{1}{2\tn(\tn-1)}$ if $i \notin \{i^*,j^*\}$, $P^{u}_{ij} = \frac{1}{2\tn^2}$ if $i\in \{i^*,j^* \}$.
Calculating $P^r$ requires us to find effective resistances on the graph. Following definition of resistance allows us to calculate them using Cayley's formula for complete graphs, 
\begin{equation*}
    R_{ij}= \frac{\text{ \# of spanning trees passing  through $(i,j)$}}{\text{ \# of spanning trees}}.
\end{equation*}
\indent 
A complete graph with $\tilde{n}$ vertices has $\tn^{\tn-2}$ spanning trees, therefore barbell graph has $\tn^{2\tn-4}(\tn^{\tn-2}\times \tn^{\tn-2})$ spanning trees. Let $K$ be the number of trees passing from an edge then  $K \times \binom{\tn}{2}=\tn^{\tn-2}(\tn-1) $. So we have $K=2\tn^{\tn-3}$. This implies that number of spanning trees passing from an edge is $2\tn^{2\tn-5}$ on barbell graph, and definitely the number of spanning trees passing from the edge $(i^*,j^*)$ is $\tn^{2\tn-4}$. This implies, $R_{ij}=1$ if $(i,j) \in \{(i^*,j^*),(j^*,i^*) \} $, $R_{ij}=\frac{2}{\tn}$ otherwise. 
Once we have explicit characterizations of $P^u$ and $P^r$, using Lemma \ref{lemma:WP-formula} we can compute the entries of $\overline{W}_{P^u}$ and $\overline{W}_{P^r}$ to be given as in $(i)$ and $(ii)$. The second largest eigenvalues of  
$\bar{W}_{P^u}$ and $\bar{W}_{P^r}$ follow from Corollary \ref{cor: stochastic e-values}.
\end{proof}

\begin{lemma}\label{lem-hitting-ub}\cite[Eqn. (2.2)]{ikeda2009hitting} Let $W$ be the transition matrix of a Markov chain with stationary distribution $\pi$. Let $j$ be a neighbor of $i$, i.e. $j\in \cN_i$, then
	$ \hij \leq  (\pi_j W_{ji})^{-1}.$
\end{lemma}

\bibliographystyle{unsrt}
\bibliography{refs_Mert,refs_Aybat,refs_Bugra}

\onecolumn

\begin{center}
     {\Huge Supplementary File}
\end{center}

\section*{\large \mg{Discussions on The Momentum-Based Acceleration Methods and ER-based Gossiping}}

\mg{ In the literature, there have been two main approaches to improve the performance of gossiping algorithms: (i) improving the communication weights, (ii) 
\sa{modifying} \mg{the} averaging scheme, \sa{e.g., adding a momentum term}. \mg{ER-based approach corresponds to the first category whereas} the papers \cite{loizou2018accelerated,loizouRevisiting,loizouAccConsensus} belong to the second category and proposes alternative averaging techniques based on a momentum term. In momentum-based approaches, the next iterate $y_{i}^{k+1}$ at node \sa{$i\in\mathcal{N}$ does not only depend on the current iterate $y_i^{k}$ but also on the previous iterate $y_i^{k-1}$ as well as $\{y_j^{k},~y_j^{k-1}\}_{j\in\mathcal{N}_i}$, i.e., the current and previous iterates of the neighbors of node $i$},} \mg{(see for example \cite{loizou2018accelerated}).}\\[0.5mm]

\bc{\bc{In the following discussion, we illustrate the benefits of momentum-based approaches and how they can be used together with effective resistance \sa{weights} to improve performance. \sa{For the sake of simplicity of the argument, we} consider the case when the updates are synchronous. In this case, if $y_i^k\sa{\in\reals}$ denotes the local estimate of the global average, $\frac{1}{n}\mathbf{1}^\top y^0$, at node $i$ in iteration $k\geq 0$, where $\mathbf{1}$ denotes the vector of ones}, gossiping algorithms consist of updates of the form:
\begin{equation}
    y^{k+1} = W y^k, \quad y^k = \begin{bmatrix} y_1^k & y_2^k & \sa{\ldots} & y_n^k \end{bmatrix}^{\sa{\top}},
    \label{eq-consensus-deter-update}
\end{equation}
starting from the initial point $y^0\sa{\in\reals^n}$, where $W$ is a doubly stochastic matrix. A common choice \sa{for the mixing matrix $W$} is 
\begin{equation} 
W = I - \alpha L 
\label{def-W-matrix}
\end{equation} where \sa{$L=[L_{ij}]_{i,j\in\mathcal{N}}$} is a symmetric weighted Laplacian matrix 
and $\alpha>0$ is a scalar satisfying $\alpha < 2/\|L\|$ (see e.g. \cite[Section 2.4]{shi2015extra}). \sa{For each $i\in\mathcal{N}$,} $L_{ij}<0$ \sa{for all $j\in \mathcal{N}_i$,}
where $\mathcal{N}_i$ is the set of neighbors of the node $i\in\mathcal{N}$; $L_{ij}=0$ if $j \not\in \mathcal{N}_i\sa{\cup\{i\}}$ and $L_{ii} = -\sum_{j\sa{\in\mathcal{N}_i}} L_{ij}>0$. Different choices of the matrix $L$ gives different algorithms. For example, \sa{\it uniform gossiping} 
corresponds to the choice \sa{$L=L^u\in\reals^{n\times n}$ such that}\footnote{\sa{Note that instead of  $L_{ij} = \frac{1}{d_i}$, for uniform gossip we set it as in \eqref{eq-laplacian-uniform} so that $L$ becomes symmetric.}}
\begin{equation} L_{ij}^u := - \frac{1}{2} \left( \frac{1}{d_i} + \frac{1}{d_j} \right), \quad  \sa{\forall j \in \mathcal{N}_i.} 
\label{eq-laplacian-uniform}
\end{equation}
Similarly, we can \sa{study gossiping based on the ER-based weights in synchronous setting by considering the choice $L=L^r\in\reals^{n\times n}$ such that}
\begin{equation} L_{ij}^r =- \frac{1}{2}\left(\frac{R_{ij}}{R_i}+\frac{R_{ij}}{R_j}\right)\quad \sa{\forall j \in \mathcal{N}_i,}
\label{eq-laplacian-er}
\end{equation}
where $R_{ij}$ is the effective resistance on the edge $(i,j)\sa{\in\mathcal{E}}$, and $R_i:=\sum_{j\in\mathcal{N}_i}R_{ij}$ \sa{for all $i\in\mathcal{N}$}.\\[0.5mm] }
\bc{
Gossiping algorithms with weighted Laplacian matrix $L$ are related to \sa{first-order, i.e., gradient-based,} optimization algorithms. To illustrate this point further, consider the \sa{following convex} quadratic optimization problem:
\begin{equation}
\min_{y\in\mathbb{R}^n}\quad \frac{1}{2}(y-y_*)^{\top} L(y-y_*),
\label{eq-quad-opt}
\end{equation}
where $y_* := \bar{y} \textbf{1}\sa{\in\reals^n}$ and $\bar{y} = \frac{1}{n}\sum_{i=1}^n y_i^0\sa{\in\reals}$ \sa{is the global average that we want to compute}. Noting that $Ly_* = \bar{y} L \textbf{1} = 0$, the updates \sa{given in} \eqref{eq-consensus-deter-update} with the choice of $W$ \sa{as} in \eqref{def-W-matrix} can be viewed as \sa{applying} a gradient descent update \sa{with step size $\alpha>0$} on the quadratic optimization problem \sa{in \eqref{eq-quad-opt}.}
From the standard theory of gradient descent, it is well-known that the distance of $y^k$ converges to $y_*$ linearly at a rate $\rho(L):=1-\alpha \lambda_{\min}^+(L)\sa{\in [0,1)}$ for $\alpha < 2/\|L\|$ where $\lambda_{\min}^+(L)$ denotes the minimum positive eigenvalue of $L$, i.e., \sa{the second smallest eigenvalue for connected graphs}. \sa{Therefore, we get the following non-asymptotic convergence:} 
$$
\Vert y^{k}-\bar{y}\sa{\mathbf{1}}\Vert^2 \leq \left(\rho(L)\right)^{2k}\Vert y^{0}-\bar{y}\sa{\mathbf{1}}\Vert^2 \quad \sa{\forall k\geq 0}. 
$$
If we set the stepsize as $\alpha = \frac{1}{\lambda_{\max}(L)}$ 
where $\lambda_{\max}(L)=\|L\|$ denotes the largest eigenvalue of $L$, we get 
$$\rho(L) = 1-\frac{1}{\kappa(L)}, \quad \mbox{where} \quad \kappa(L) := \sa{\frac{\lambda_{\max}(L)}{\lambda_{\min}^+(L)}}$$
is called the condition number. When the condition number $\kappa(L)$ is very large, the convergence can be slow.
Adding a momentum term is a technique to improve the convergence rate of gradient descent methods with respect to its dependency to the condition number. For example, Polyak's heavy-ball (HB) method applied to the objective \eqref{eq-quad-opt} consists of the iterations
 \begin{equation} \label{eq-heavy-ball}
 y^{k+1} = (I - \alpha L)y^k + \beta(y^k - y^{k-1}),
\end{equation}
where the last term $\beta(y^k - y^{k-1})$ is referred to as the \emph{momentum term} and $\beta$ is called the momentum parameter (see e.g., \cite{loizou2018accelerated}). The convergence rate of the heavy-ball (HB) method on quadratic objectives of the form \eqref{eq-quad-opt} has been well-studied in the literature and it can be shown that the heavy-ball method given in iterations \eqref{eq-heavy-ball} will converge to the consensus vector $y_*$ with the asymptotic linear convergence rate 
 $$
 \rho_{HB}(L):=1- \Theta \left(\frac{1}{\sqrt{\kappa(L)}} \right)
 $$
for a specific choice of the stepsize $\alpha$ provided that $\beta$ is tuned properly as a function of the eigenvalue \sa{$\lambda^+_{\min}(L)$} \cite{loizou2018accelerated}. Achieving this rate with the choice of $\beta$ in \cite{loizou2018accelerated} would require estimating $\lambda_{\min}^+(L)$. That being said for ill-conditioned problems when the condition number $\kappa(L)$ is sufficiently large, we observe that HB converges faster, i.e. $\rho_{HB}(L)< \rho(L)$. For example, for the barbell graph, with an analysis similar to that in Proposition 8 of the revised manuscript, we can characterize the eigenvalues of the weighted graph Laplacians $L^u$ and $L^r$ that correspond to uniform weights and ER-based weights given in \eqref{eq-laplacian-uniform} and \eqref{eq-laplacian-er} respectively and obtain
\begin{equation} 
\kappa(L^u) = \Theta({n^2}), \quad \kappa(L^r) = \Theta(n).
\label{eq-spec-gap-er-vs-uniform}
\end{equation}
Therefore, without momentum averaging (when $\beta=0$), we obtain the convergence rates
\begin{align}
    \rho^{u}:= \rho(L^u) = 1-\Theta(\frac{1}{n^2}),\;\; \rho^{r}:= \rho(L^r) = 1-\Theta(\frac{1}{n}),
    \label{eq-barbell-rate-comparison}
\end{align}
for uniform weights and ER-based weights. On the other hand, for HB method, we obtain the rates 
\begin{align*}
    \rho_{HB}^{u}:= \rho_{HB}(L^u) = 1-\Theta(\frac{1}{n}),\;\; \rho_{HB}^{r}:= \rho_{HB}(L^r) = 1-\Theta(\frac{1}{\sqrt{n}}).
\end{align*}
We observe that fastest rate is obtained by using the HB method on the quadratic problem in \eqref{eq-quad-opt} defined by the weighted Laplacian corresponding to the ER-weights, i.e., $\rho_{HB}^r$ is the fastest rate in terms of its dependency to $n$. This shows that ER weights can be used together with momentum averaging techniques. Basically, from  \eqref{eq-spec-gap-er-vs-uniform}, we observe that effective-resistance based approach yields to a better conditioned Laplacian compared to uniform weights; and further improvement can be achieved by employing momentum averaging. In other words, ER weights are needed to improve the conditioning of the weighted Laplacian matrix and momentum-based approaches can be used on top of this to get further performance improvement. Besides the HB method, Nesterov's accelerated gradient method is an alternative momentum averaging-based technique which will also yield to similar accelerated convergence rates.\\[0.5mm]
}
\bc{
The discussion we provided was for the synchronous setup, the asynchronous setup can be analyzed similarly.\footnote{\mg{In the asynchronous setup, at every iteration, node $i$ contacts a neighbor randomly to update its decision variable rather than contacting all the neighbors. In the case of the barbell graph, each node has $\Theta(n)$ neighbors so needs on average $\Theta(n)$ iterations to contact all the neighbors. Consequently, more iterations will be required to converge compared to the synchronous setup. With a similar analysis to above, it can be shown that ER weights \bc{on barbell graphs lead} to
$\mathbb{E}\|y^k - y^*\|^2 \leq \bc{\left[\rho_{async}^r\right]^{2k}} \mathbb{E} \|y^0 - y^*\|^2$
where $\rho^r_{async} = 1 - \Theta(\frac{1}{n^3})$ instead of $\rho^r = 1 - \Theta(\frac{1}{n^2})$ obtained above in \eqref{eq-barbell-rate-comparison}. The rate $\rho^r_{async}$ also follows directly from \bc{Proposition 9}.} }
}

\section*{\large \mg{Further Discussions on Our Conductance Bounds and Averaging Time with Effective Resistances}}

\mg{We recall that the averaging time $T_{ave}(\varepsilon,P)$ with an 
expected iteration matrix $\bar{W}_{P}$ satisfies \begin{equation} 
T_{ave}(\varepsilon,P)=\Theta\left(T(\overline{W}_P) \right),\quad \mbox{where}\quad T(\overline{W}_P):=  \frac{1}{\log (1/\lambda_{n-1}(\overline{W}_P))}, 
\label{eq-time-averaging}
\end{equation} where $\lambda_{n-1}(\cdot)$ denotes the second-largest eigenvalue. 
Therefore, comparing effective-resistance (ER) weights with uniform weights amounts to comparing the second-largest eigenvalues $\lambda_{n-1}(\overline{W}_{P^r})$ and $\lambda_{n-1}(\overline{W}_{P^{u}})$, where \sa{$\overline{W}_{P^r}$ and $\overline{W}_{P^u}$ are the expected iteration
matrices defined using ER and uniform weights, respectively.} For barbell graphs (that correspond to the special case of \sa{$c$-barbell graphs} with $c=2$), our analysis is tight as we have developed an explicit formula for computing the second-largest eigenvalue of the matrix $\overline{W}_{P^r}$ as well as the second-largest eigenvalue of  $\overline{W}_{P^u}$}.
\bc{However, for $c$-barbell graphs with $c>2$, the second-largest eigenvalues of the gossiping matrices $\overline{W}_{P^r}$ and $\overline{W}_{P^u}$ are not explicitly known.
Therefore, in our paper, we resorted to the conductance bounds which is a common technique in the literature to obtain lower and upper bounds on the second \sa{largest} eigenvalue $\lambda_{n-1}(\overline{W}_{P})$ and consequently the spectral gap $ \Delta := 1-\lambda_{n-1}(\overline{W}_{P})$ through the Cheeger inequalities. Based on this approach, we can obtain the following lower and the upper bounds for the spectral gaps $\Delta_{r} \sa{:=} 1 - \lambda_{n-1}(\overline{W}_{P^r})$ and $\Delta_{u}  := 1 - \lambda_{n-1}(\overline{W}_{P^u})$ that correspond to ER and uniform weights, respectively:
\begin{align*}
{\Delta}^{{lb}}_{r} &:= \Phi^2(\bar{W}_{P^r}) \leq {\Delta}_{r} \leq {\Delta}^{{ub}}_{r}:= 2\Phi(\bar{W}_{P^r}),\\
{\Delta}^{{lb}}_{u} &:= \Phi^2(\bar{W}_{P^u}) \leq {\Delta}_{u}  \leq {\Delta}^{{ub}}_{u} := 2\Phi(\bar{W}_{P^u}),
\end{align*}
}
\sa{where $\Phi(\overline{W}_{P})$ denotes the graph conductance as defined in the paper for the reversible Markov chain corresponding to transition probability matrix $\overline{W}_{P}$.}

\bc{To illustrate the tightness of our bounds, we consider \emph{the approximation ratio}, i.e., the ratio of these bounds in a logarithmic scale
     $${a}_{r}^{lb}:= \frac{\log (\Delta^{lb}_{r})}{\log (\Delta_{r})} \leq 1 \leq 
     a^{ub}_{r}:= \frac{\log ({\Delta}^{{ub}}_{r})}{\log ({\Delta_{r}})}.
     $$     
We define $a_u^{lb}$ and $a_u^{up}$ similarly for the uniform weights.}
     
\begin{table}[h!]
\centering
\begin{tabular}{|c|c|c|c|c|}
\hline
\textbf{$\tilde{n}$} & \textbf{${a}_{r}^{lb}$} & \textbf{${a}_{r}^{ub}$} & \textbf{${a}_{u}^{lb}$} & \textbf{${a}_{u}^{ub}$} \\ \hline
10 & 0.802 & 1.735 & 0.865 & 1.848 \\ \hline
16 & 0.818 & 1.756 & 0.883 & 1.873 \\ \hline
18 & 0.822 & 1.761 & 0.887 & 1.878 \\ \hline
20 & 0.825 & 1.765 & 0.891 & 1.882 \\ \hline
22 & 0.828 & 1.769 & 0.893 & 1.886 \\ \hline
28 & 0.834 & 1.778 & 0.900 & 1.894 \\ \hline
30 & 0.836 & 1.780 & 0.901 & 1.896 \\ \hline
36 & 0.841 & 1.786 & 0.905 & 1.900 \\ \hline
38 & 0.842 & 1.788 & 0.906 & 1.902 \\ \hline
44 & 0.845 & 1.793 & 0.909 & 1.905 \\ \hline
46 & 0.846 & 1.794 & 0.910 & 1.906 \\ \hline
48 & 0.847 & 1.795 & 0.911 & 1.907 \\ \hline
50 & 0.848 & 1.797 & 0.912 & 1.908 \\ \hline
100 & 0.862 & 1.815 & 0.923 & 1.920 \\ \hline
500 & 0.886 & 1.847 & 0.939 & 1.938 \\ \hline 
1000& 0.894 & 1.858 &0.944 &  1.943 \\ \hline
\end{tabular}
\caption{\bc{Comparison of $\log(2\Phi(\overline{W}_P))/\log(1-\lambda_{n-1}(\overline{W}_P))$ and $\log(\Phi^2(\overline{W}_p)))/\log(1-\lambda_{n-1}(\overline{W}_{P}))$ of ER based gossiping and classical gossiping on the ($c$-barbell) $c-K_{\tilde{n}}$ graph with $c=10$.}}
\label{fig:Cond_Comparison}
\end{table}
\bc{
The closer the ratios ${a}_{r}^{lb}$ and ${a}_{r}^{ub}$ are to 1, the better the approximation quality is. In Table \ref{fig:Cond_Comparison}, we illustrate the tightness of our bounds for $c$-barbell graphs $(c-K_{\tilde n})$ that consists of $c$ cliques where each clique has $\tilde{n}$ nodes,  where we report ${a}_{r}^{lb}$, ${a}_{r}^{ub}$. We also display the ratios ${a}_{u}^{lb}$ and ${a}_{u}^{ub}$ for uniform weights, which are computed similarly. The results illustrate that all the ratios lie in a reasonable range (in the interval $[0.80,1.95]$) with lower bounds being tighter than the upper bounds.
These results show that conductance-based analysis leads to useful approximations. In particular, we can see that the lower bounds are becoming tighter ($a_r^{lb}$ is increasing) as the number of nodes increases on the graph.
}
\mg{As an additional experiment, we also computed the eigenvalues of $\overline{W}_{P^r}$ and $\overline{W}_{P^u}$ with the standard eigenvalue solver in Matlab 2021a (using the function \emph{eig} with default settings).
\sa{Using the second largest eigenvalues of $\overline{W}_{P^r}$ and $\overline{W}_{P^u}$}, we compute the times $T(\overline{W}_{P^r})$ and $T(\overline{W}_{P^{u}})$ required for both approaches. From \eqref{eq-time-averaging}, we see that 
\begin{equation} \sa{\frac{T_{ave}(\epsilon,P^u)}{T_{ave}(\epsilon,P^r)}=\Theta\left(\frac{T(\overline{W}_{P^{u}})}{T(\overline{W}_{P^{r}})}\right),}\quad \frac{T(\overline{W}_{P^{u}})}{T(\overline{W}_{P^{r}})} =  
\frac{\log ([\lambda_{n-1}(\overline{W}_{P^r})])}{\log ([\lambda_{n-1}(\overline{W}_{P^u})])}.
\label{eq-ratio-times}
\end{equation}
In Figure \ref{fig:spec_cap_c_barbell}, we plot the ratio on the right hand-side for the $c$-barbell graph, denoted as \bc{$c-K_{\tilde{n}}$}. For different values of $c$ fixed, we vary $\tilde{n}$ and observe that the ratio $\frac{\log ([\lambda_{n-1}(\overline{W}_{P^r})])}{\log ([\lambda_{n-1}(\overline{W}_{P^u})])}$ is always \bc{larger} than 1 and the ratio is growing as $\tilde{n}$ increases. This shows that ER weights admits better (smaller) averaging times for especially large networks, i.e., the performance gain being more and more significant as the number of nodes $\tilde{n}$ increases. In light of these experiments, we can conclude the superiority of ER weights \sa{over the uniform weights} from a numerical perspective as well.
\begin{figure}[h!]
    \centering
    \includegraphics[width=0.7\linewidth]{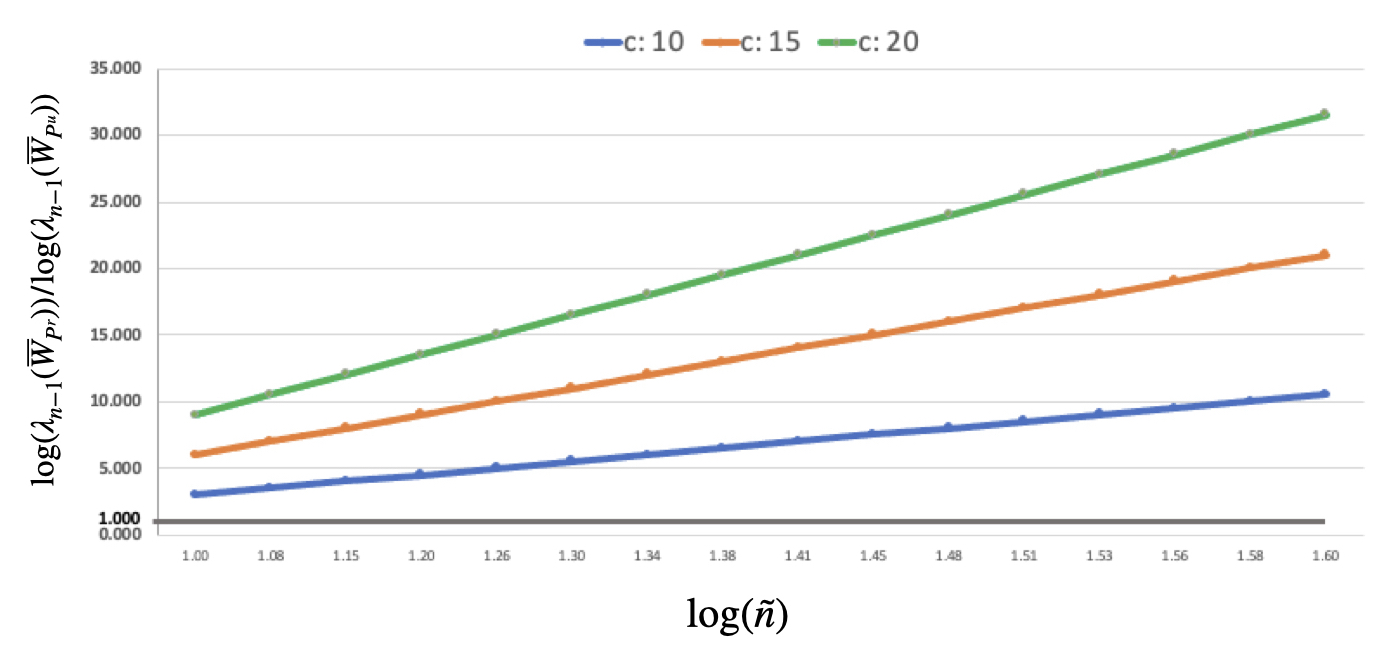}
    \caption{ \bc{The comparison of the ratio between the times $T(\overline{W}_{P^u})$ and $T(\overline{W}_{P^r})$ on $c-K_{\tilde{n}}$.} }
    \label{fig:spec_cap_c_barbell}
\end{figure}     
}

\section*{\large Normalized D-RK Algorithm}\label{sec-method} 
D-RK method for computing the effective resistances in a decentralized way and its normalized version which we call normalized D-RK has been introduced in \cite{aybat2017decentralized} where the authors show that these methods converge linearly with rates
\begin{small}
\BEQ\label{def-rho}
\rho \triangleq 1-\left(\frac{\lambda_{\min}^+(\cL)}{\norm{\cL}_F}\right)^2, \quad {\rho_S} \triangleq 1-\frac{1}{n}\lambda_{\min}^{+}(\cL S^{-1}\cL),
\EEQ
\end{small}
respectively where $\lambda_{\min}^+(\cdot)$ denotes the smallest positive eigenvalue and $S$ is a normalization matrix defined as 
\begin{small}
\BEQ\label{def-S} S=\diag(s) \quad \mbox{where} \quad s_i\triangleq\sum_{j\in\cN_i\cup\{i\}}\cL_{ij}^2 \mbox{ for } i\in\cN.
\EEQ
\end{small}
{Based on numerical evidence, it was conjectured in \cite{aybat2017decentralized} that normalized D-RK is faster than D-RK, i.e. $\rho_S\leq \rho$.}
First, we provide a technical result and then the following proposition proves this conjecture.
\begin{lemma}\label{lemma-si}The Laplacian $\mathcal{L}$ has the following property:
$\frac{1}{n^2}\sum_{i=1}^{n}\frac{1}{s_i}\geq \frac{1}{||\mathcal{L}||_{F}^{2}}$,
where $s_i$ is defined by \eqref{def-S}.
\end{lemma}
\begin{proof} 
Note that $||\mathcal{L}||_{F}^{2}=\sum_{i=1}^{n}\sum_{j=1}^{n}\mathcal{L}_{ij}^{2}=\sum_{i=1}^{n}\sum_{j\in N_{i}\cup\{i\}}\mathcal{L}_{ij}^2=\sum_{i=1}^{n}s_{i},$ where we used the fact that $\mathcal{L}_{ij}=0$ for all $(i,j)\notin \cE$. Applying arithmetic-harmonic mean inequality to the sequence $\{s_{i}\}_{i \in \{ 1,..,n \}}$, we obtain
 $	\frac{1}{n}||\mathcal{L}||_{F}^{2}=\frac{1}{n}\sum_{i=1}^{n}s_i \geq n \Big[\sum_{i=1}^{n}\frac{1}{s_i}\Big]^{-1}$.
We conclude by multiplying both sides with $1/n$.
\end{proof}
Now we are ready to prove our conjecture.
\begin{proposition}\label{prop-rate-compare} For $S$ defined by \eqref{def-S}, the following inequality holds:
$
\frac{1}{n}\lambda_{\min}^+ (\cL S^{-1} \cL) \geq \left(\frac{\lambda_{\min}^+ (\cL)}{\| \cL\|_F}\right)^2$. Then, it follows that $\rho_S \leq \rho$ where $\rho$ and $\rho_S$ are defined by \eqref{def-rho}. 
\end{proposition}
\begin{proof} 
Since $\mathcal{L}$ and $S$ are symmetric matrices so are $\mathcal{L}^{2}$ and $S^{-1}$. Let $\{\lambda_i(\mathcal{L})\}_{i=1}^n$ and $\{\lambda_i(\mathcal{S})\}_{i=1}^n$ denote the eigenvalues of these matrices sorted in increasing order, i.e. $\lambda_n$ is the largest eigenvalue, $\lambda_1$ is the smallest one. By the eigenvalue interlacing result in \cite[Chapter 2, Eq. (2.0.7)]{Zhang}, we obtain\footnote{We set $l=n$ and $i_t=2$ for $t=1,\ldots,l$ in Eq.~(2.0.7) in~\cite{Zhang}.}
\begin{equation}\label{eq: specialinterlace}
n \lambda_{2}(\mathcal{L}^2\mathcal{S}^{-1}) \geq\lambda_{2}(\mathcal{L}^{2})\sum_{i=1}^{n}\lambda_{i}(\mathcal{S}^{-1}),	
\end{equation}
where all the matrices have non-negative real eigenvalues as both $\mathcal{L}$ and $S$ are symmetric with non-negative eigenvalues. 
Clearly, $\lambda_{2}(\mathcal{L}^{2})=\lambda_{2}(\mathcal{L})^{2}>\lambda_{1}(\mathcal{L}^2)=0$. 
Furthermore, the eigenvalues of $\mathcal{L}^2 \mathcal{S}^{-1}$ and $\mathcal{L} \mathcal{S}^{-1}\mathcal{L}$ are the same \footnote{If $u$ is an eigenvector of the latter matrix corresponding to a non-zero eigenvalue $\lambda$, then $\mathcal{L}u$ would be the right eigenvector of the former matrix with the same eigenvalue; similarly, if $u$ is a right-eigenvector of $\mathcal{L}^2 \mathcal{S}^{-1}$ corresponding to a nonzero eigenvalue $\lambda$, then $\cL \cS^{-1} u$ is an eigenvector of $\mathcal{L} \mathcal{S}^{-1}\mathcal{L}$ with the same eigenvalue.} Therefore, since $\mathcal{L} \mathcal{S}^{-1}\mathcal{L}$ is positive semidefinite with $\lambda_1(\mathcal{L} \mathcal{S}^{-1}\mathcal{L})=0$, we also have
\begin{equation}\label{eq-zero-eig}
\lambda_1(\mathcal{L}^2 \mathcal{S}^{-1})=0.
\end{equation}%
Moreover, $S$ is a diagonal matrix with diagonal entries $S_{ii}=s_i$; therefore, eigenvalues of $S$ are given by $s_i$ with $i=1, 2, \dots, n$. Hence \eqref{eq: specialinterlace} is equivalent to 
\begin{eqnarray}
n \lambda_{2}(\mathcal{L}^2\mathcal{S}^{-1})  \geq\lambda_{2}(\mathcal{L})^{2}\sum_{i=1}^{n}\frac{1}{s_i}
\geq \lambda_{\min}^{+}(\mathcal{L})^{2} \frac{n^2}{\|\mathcal{L}\|_F^2}> 0, \label{ineq-eigenvalue-S-b} 
\end{eqnarray}
where the inequalities follow from Lemma \ref{lemma-si} and the fact that  $\lambda_{2}(\mathcal{L}) =\lambda_{\min}^{+} (\mathcal{L})>0$ 
due to $\cG$ being a connected graph, where $\lambda_{\min}^{+}(\cdot)$ denotes the smallest positive eigenvalue. From \eqref{eq-zero-eig} and \eqref{ineq-eigenvalue-S-b}, we conclude that $\lambda_{2}(\mathcal{L}^2\mathcal{S}^{-1})$ is the smallest positive eigenvalue of $\mathcal{L}^2\mathcal{S}^{-1}$, i.e., 
\BEQ\label{eq-lambda-S}\lambda_{2}(\mathcal{L}^2\mathcal{S}^{-1}) = \lambda_{\min}^{+} (\mathcal{L}^2\mathcal{S}^{-1}).\EEQ   
Finally, 
using the fact that the eigenvalues of $\mathcal{L}^2 \mathcal{S}^{-1}$ and $\mathcal{L} \mathcal{S}^{-1}\mathcal{L}$ are the same once again, we get 
$\lambda_{\min}^{+}(\mathcal{L}\mathcal{S}^{-1}\mathcal{L}) = \lambda_{\min}^{+}(\mathcal{L}^2\mathcal{S}^{-1})$.
Combining this with \eqref{ineq-eigenvalue-S-b} and \eqref{eq-lambda-S} leads to 
\begin{equation*}
\frac{1}{n}\lambda_{\min}^{+}(\mathcal{L}\mathcal{S}^{-1}\mathcal{L})= \frac{1}{n}\lambda_{\min}^{+}(\mathcal{L}^2\mathcal{S}^{-1})\geq \Big( \frac{\lambda_{\min}^{+}(\mathcal{L})}{||\mathcal{L}||_{F}} \Big)^{2},
\end{equation*}
which directly implies $\rho_S\leq \rho$. This completes the proof. 
\end{proof} 

\section*{\large Proof of Proposition~\ref{prop: gen. eigenvalue}}

 The proof follows by adapting the proof of \cite[Proposition 5.1]{BarbellBoyd} to our setting with minor modifications. It is based on exploiting the symmetry group properties of the barbell graph with algebraic techniques. We first give relevant background material below before going into the details of the proof. 
{
\subsection*{Background Material}
Consider a weighted graph $\cG=(\cN,\cE,w)$. 
A \emph{permutation} $p: \mathcal{N}\rightarrow \mathcal{N}$ is a mapping that rearranges the vertices, i.e. it is a bijection from the node set $\cN$ to itself. We consider a \emph{permutation group} $H$, which is a group whose elements are permutations of $\cN$ and whose group operation is the composition of permutations in $H$. By the group property, if two permutations $s_1, s_2 \in H$, then the composition $s_1 s_2 \in H$ and in particular the \emph{identity permutation} $e$ which maps all the elements of $\cN$ to itself is also contained in $H$. The group that contains all the $n!$ permutations defined on $\cN$ is denoted as $S_{n}$.

The \emph{direct product} $(H_1\times H_2)$ of two groups $H_1,H_2$ is defined as the group that consists of elements from the Cartesian product of $H_1$ and $H_2$ with the elementwise composition, i.e. $(h_1, h_2) \in (H_1\times H_2)$ if and only if $h_1 \in H_1$ and $h_2 \in H_2$ and if $(h_1,h_2) \in (H_1 \times H_2) $ and $(\tilde{h}_1,\tilde{h}_2) \in (\tilde{H}_1 \times \tilde{H}_2) $ then the composition operation $\cdot$ over $(H_1\times H_2)$ is defined as $(h_1,h_2) \cdot (\tilde{h}_1,\tilde{h}_2)=(h_1\tilde{h}_1,h_2\tilde{h}_2)$.  A subgroup $M$ of a group $H$ is \textit{normal} if for all $h \in H$ and $m \in M$ we have $hmh^{-1}\in M$. The \textit{semidirect product} $H_1\ltimes H_2$ of two groups $H_1$ and $H_2$ is the group that consists of elements $h=h_1h_2$ with $h_1 \in H_1$ and $h_2 \in H_2$ and the subgroup $H_1 $ 
is normal in $H_1\ltimes H_2$ with the condition $H_1\cap H_2=\{e\}$. The \emph{orbit} $O_{i}$ of an element $i \in \mathcal{N}$, under a permutation group $H$ is  the set $O_{i}\triangleq \{v\in \cN ~|~ \exists s \in H \;\;\text{s.t.}\;\; s(v)=i\}$. In other words, the orbit of node $i$ is the set of vertices that can be mapped to $i$ by an element of the permutation group $H$. This definition creates an equivalence relation $\sim$ on $\cN$; for $i,j \in \cN$, we say $i \sim j$ if $O_i = O_j$. In particular, equivalence classes form a partition of $\cN$. 

A permutation $s$ is called \textit{an automorphism} of the weighted graph $\cG$ if the weight matrix $w$ 
is invariant under $s$, i.e. if $w(i,j)=w(s(i),s(j))$. From this definition, an automorphism $s$ also satisfies $W(i,j)=W(s(i),s(j))$ where $W(i,j) = w(i,j)/\sum_{j\in \cN_i} w(i,j)$ is the transition probability. We are interested in such permutations that preserve the structure of $w$ and therefore $W$. The group of all automorphisms with the operation of composition of permutations is called the \emph{automorphism group} of the graph and is denoted by  $\mbox{Aut}(\cG)$. Let $S$ be a subgroup of $\mbox{Aut}(\cG)$ and consider the orbits $\{O_i\}_{i\in\cN}$ under the permutation group $S$ which partition the set $\cN$.
We define \emph{orbit graph} to be the graph whose vertices consist of the equivalence classes $O_i$ for $i\in\cN$ and we consider an induced Markov chain on the orbit graph with probability transition probabilities defined as
\begin{equation}
\label{def: inducedweights}
W_{S}(O_i,O_j)=\sum_{j'\in O_{j}}W(i,j'). 
\end{equation}
This Markov chain is also called the \emph{orbit chain}. It can be shown that the definition of the weights $W_S$ above does not depend on the choice of the element $i$ from the set $O_i$ (see e.g. \cite{BarbellBoyd}). 
}
\\
\subsection*{Proof}
First, we consider the automorphism group of the barbell graph $K_{\tn}-K_{\tn}$ with edge weights given by Proposition \ref{prop: gen. eigenvalue}. Consider the nodes $i_*$ and $j_*$ that connect the complete subgraphs of the barbell graph and without loss of generality assume that we enumerate the nodes so that $i_*= \tn$, $j_* = \tn+1$ and a node $i<\tn$ is on the complete subgraph on the left hand-side and any node $j>\tn+1$ is on the complete subgraph on the right-hand side. We see from the symmetry structure of $W$ that if we take any two nodes from a complete subgraph and permute them, this would be an automorphism. Similarly, swapping the two complete subgraphs between them would be an automorphism; i.e. the permutation  $C_2:\mathcal{N}\rightarrow \mathcal{N}$ that maps $i \overset{C_2}{\mapsto}-i~\mod(n+1)$ is an automorphism. It follows from these observations that
the automorphism group of $K_{\tn}-K_{\tn}$ is the group $C_2 \ltimes (S_{\tn-1}\times S_{\tn-1})$ (see also \cite{BarbellBoyd} for more details). It is known that for any subgroup $S$ of the automorphism group, the eigenvalues of the transition matrix $W_S$ defined by \eqref{def: inducedweights} should also be an eigenvalue of the transition matrix $W$ (see e.g. \cite[Section 3]{BarbellBoyd}). Note that the square matrix $W_S$ has dimension $n_S \times n_S$ where $n_S \leq n$, so the set of eigenvalues of $W_S$ are a subset of the set of all eigenvalues of $W$. 
We are going to use this result to prove the Proposition \ref{prop: gen. eigenvalue}. Next, we consider the eigenvalues of the transition matrices $W_S$ of the orbit chains under subgroups $S$ of $C_2\ltimes ( S_{\tn-1}\times S_{\tn-1})$:
\begin{figure}[ht!]
\centering
\includegraphics[width=0.30\linewidth]{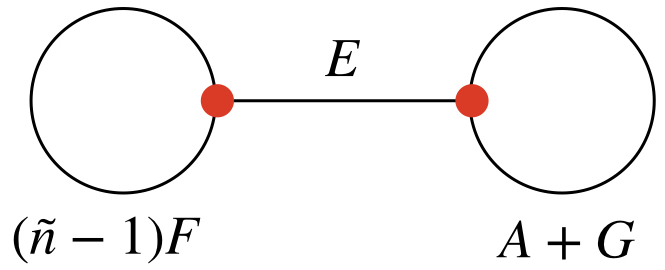}
\caption{ Orbit graph under $C_2 \ltimes (S_{\tilde{n}-1}\times S_{\tilde{n}-1})$}
\label{fig: Orbitchaina}
\end{figure}

a) The orbit chain under $C_2 \ltimes (S_{\tn-1}\times S_{\tn-1})$ (Figure \ref{fig: Orbitchaina}) has the transition matrix
$
\begin{bmatrix}
\frac{A+G}{A+G+E} &\frac{E}{A+G+E}\\
\frac{E}{(n-1)F+E} & \frac{(n-1)F}{(n-1)F+E}
\end{bmatrix}$.
Since $\lambda_a=1$ is an eigenvalue, and its trace is the sum of eigenvalues; it follows that the other eigenvalue of this matrix is given by $\lambda_b= -1+\frac{A+G}{A+G+E}+\frac{F}{F+B}$.

b) Consider the orbit chain under $C_2$ illustrated on the left panel of Figure \ref{orbitsb}.
\begin{figure}[ht!]
\centering
\includegraphics[width=0.30\linewidth]{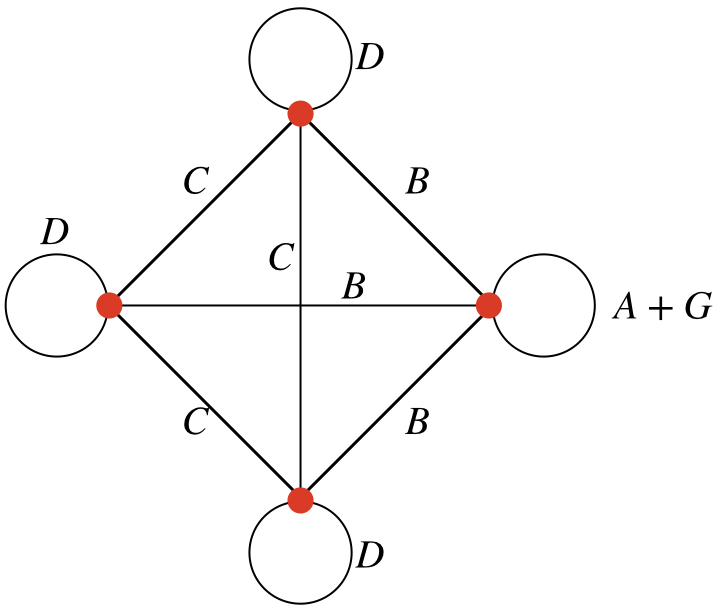}
\includegraphics[width=0.30 \linewidth]{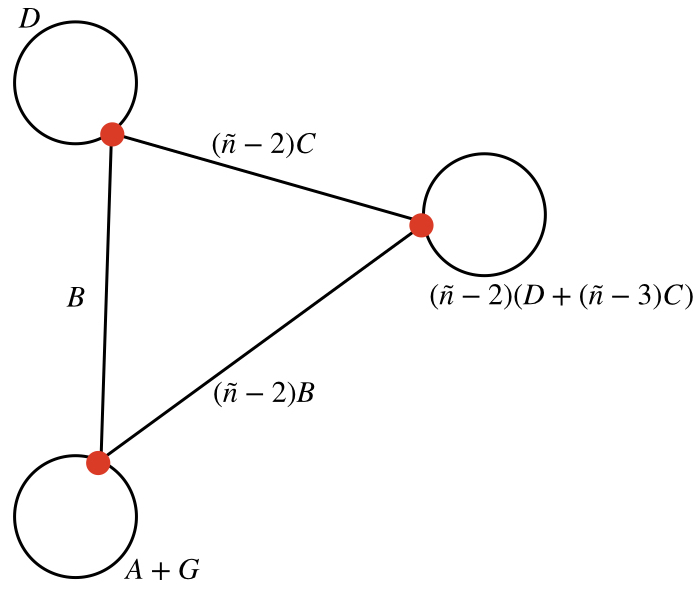}
\caption{\textbf{Left:} Orbit graph under $C_2 $. \textbf{Right:} Orbit graph under $C_2\ltimes (S_{\tn-2}\times S_{\tn-2}) $ }
\label{orbitsb}
\end{figure}
This orbit graph has two orbits under permutation  $S_{\tn-1}$: One of them contains only one node  (the node with a self-loop with weight $(A+G)$) and the other orbit has the remaining $\tn-1$ nodes. 
Notice that the latter orbit has identical $\tn-1$ elements and therefore the permutation group $C_2 \ltimes (S_{\tn-2}\times S_{\tn-2})$ fixes one of the nodes having a loop with weight $D$ and permutes the remaining $\tn-2$ nodes among themselves without affecting the orbit with one node. Therefore, by \cite[Thereom 3.1]{BarbellBoyd}, the eigenvalues of the transition matrix $W'$ of the orbit graph obtained by the permutation group $S = C_2\ltimes (S_{\tn-2}\times S_{\tn-2})$ (illustrated on the right panel of Figure \ref{orbitsb}) are also eigenvalues of the transition matrix $W$. The transition matrix $W'$ is $3\times 3$ with three eigenvalues, including $\lambda_{a}$ and $\lambda_{b}$ that  we have already found at part $(a)$.
The third eigenvalue $\lambda_c$ can be computed from the transition matrix $W'$ of the orbit chain under $C_2\ltimes (S_{\tn-2}\times S_{\tn-2})$: 
\begin{equation*}
\begin{bmatrix}
\frac{A+G}{(\tn-1)B+A+G} & * & * \\ 
* & \frac{D}{(\tn-2)C+D+B} & * \\ 
* & * & \frac{D+(\tn-3)C}{B+D+(\tn-2)C}
\end{bmatrix},
\end{equation*}
where we use $*$ to denote the entries of this matrix that will not be relevant to our discussion. In particular, the eigenvalues of this matrix will be $\lambda_a$, $\lambda_{b}$ and $\lambda_{c}$; the latter will be an eigenvalue of $W$ with multiplicity $2\tn-4$. Again, using the fact that the trace of a matrix is equal to the sum of its eigenvalues, we obtain
\[
\lambda_{c}= \frac{D-C}{F+B}.
\]
c) Lastly, orbit chain under $(S_{\tn-1}\times S_{\tn-1})$ consists of four orbits: $(\tn-1)$ points in the left and right complete graphs and vertices $i_*$ and $j_*$ as illustrated in Figure \ref{orbitsc}. 
\begin{figure}[ht!]
\centering 
\includegraphics[width=0.30\linewidth]{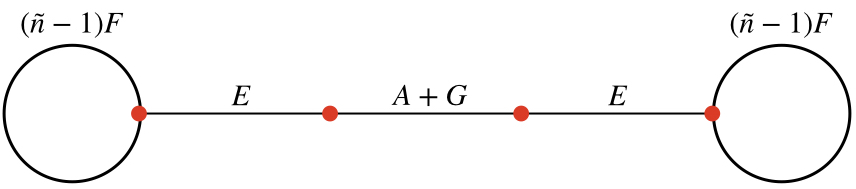}
\caption{Orbit graph under $S_{\tilde{n}-1}\times S_{\tilde{n}-1}$}
\label{orbitsc}
\end{figure}
\\
This orbit chain has the transition matrix of the form
\begin{equation*}
\begin{bmatrix}
\frac{F}{B+F}	&\frac{B}{B+F} &0 & 0\\
\frac{E}{A+E+G}	& \frac{G}{A+E+G}&\frac{A}{A+E+G} & 0 \\ 
0  &\frac{A}{A+E+G} &\frac{G}{A+E+G} & \frac{E}{A+E+G}\\ 
	0&0 &\frac{B}{B+F} &\frac{F}{B+F} 
\end{bmatrix}. 
\end{equation*}
After a straightforward computation, it can be checked that this matrix has the eigenvalues, $1, \lambda_{+},\lambda_{-}, (-1+\frac{A+G}{A+E+G}+\frac{F}{B+F})$ 
where
\begin{equation*}
\lambda_{\pm}=\frac{1}{2}\Bigg[\frac{F}{B+F}+\frac{G-A}{A+E+G} \,\pm\,\sqrt{S} \Bigg],
\end{equation*}
and $S=\bigg(\frac{F}{B+F}+\frac{G-A}{A+E+G}\bigg)^2-\frac{4(FG-BE-AF)}{(B+F)(A+E+G)}$. 
\\ 
\begin{remark} Boyd \emph{et al.} \cite{BarbellBoyd} studied the case $W_{i^*i^*}=0=W_{j^*j^*}$ where similar orbit chains and graphs arise. The proof of Proposition \ref{prop: gen. eigenvalue} given here is a minor modification of the original proof of Boyd et al. \cite[Proposition 2.2]{BarbellBoyd} and extends it to the more general case where $W_{i^*i^*}$ or $W_{j^*j^*}$ can be strictly positive.
\end{remark}


\twocolumn

\end{document}